\documentclass[10pt]{article}
\usepackage{amssymb,amsmath,amsthm,amsfonts,mathrsfs}
\usepackage[pdftex,colorlinks,linktocpage,citecolor=black,linkcolor=black]{hyperref}
\usepackage{verbatim}
\usepackage{mathtools}
\usepackage[usenames]{color}
\usepackage{graphicx,graphics}
\usepackage{epstopdf}
\usepackage{subfigure}
\usepackage{float,wrapfig}
\usepackage{bm}
\usepackage{multirow,multicol}
\usepackage{verbatim}
\usepackage{booktabs}
\usepackage{threeparttable}
\usepackage[linesnumbered,ruled,algosection]{algorithm2e}
\graphicspath{{./pics/}}
\allowdisplaybreaks

\usepackage{enumitem}
\numberwithin{equation}{section}
\newtheorem{theorem}{Theorem}[section]
\newtheorem{lemma}{Lemma}[section]
\theoremstyle{remark}
\newtheorem{remark}{Remark}[section]

\newtheorem{assumption}{Assumption}[section]

\def\cT {\mathcal T}
\def\bold {\boldsymbol}
\def\Om {\Omega}
\def\J {\mathcal{J}}
\def\F {\mathcal{F}}
\def\U {\mathcal{U}}
\def\to{\rightarrow}
\def\eps {\varepsilon}
\def\R {\mathbb{R}}

\renewcommand{\S}{\bold{\mathcal{S}}}
\newcommand{\E}{\bold{\mathcal{E}}}

\newcommand{\dx}{\,{\rm d}\bold{x}}
\newcommand{\ds}{\,{\rm d}s}
\newcommand{\dd}{\,{\rm d}}

\begin{document}

\title{\Large On the Crouzeix-Raviart Finite Element Approximation of Phase-Field Dependent Topology Optimization in Stokes Flow}

\author{Bangti Jin$^1$
\and Jing Li$^2$  \and Yifeng Xu$^3$ \and Shengfeng Zhu$^{2\,4}$}

\footnotetext[1]{Department of Mathematics, The Chinese University of Hong Kong, Shatin, New Territories, Hong Kong, China (b.jin@cuhk.edu.hk, bangti.jin@gmail.com)}
\footnotetext[2]{School of Mathematical Sciences, East China Normal University, Shanghai 200241, China. (betterljing@163.com)}
\footnotetext[3]{Department of Mathematics \& Scientific Computing Key Laboratory of Shanghai Universities, Shanghai Normal University, Shanghai 200234, China. (yfxu@shnu.edu.cn, yfxuma@aliyun.com)}
\footnotetext[4]{Key Laboratory of Mathematics and Engineering Applications \& Shanghai Key Laboratory of Pure Mathematics and Mathematical Practice, East China Normal University, Shanghai 200241, China. (sfzhu@math.ecnu.edu.cn)}

\date{}
\maketitle
\begin{abstract}
In this work, we investigate a nonconforming finite element approximation of phase-field parameterized topology optimization governed by the Stokes flow. The phase field, the velocity field and the pressure field are approximated by conforming linear finite elements, nonconforming linear finite elements (Crouzeix-Raviart elements) and piecewise constants, respectively. When compared with the standard conforming counterpart, the nonconforming FEM can provide an approximation with fewer degrees of freedom, leading to improved computational efficiency. We establish the convergence of the resulting numerical scheme in the sense that the sequences of phase-field functions and discrete velocity fields contain subsequences that converge to a minimizing pair of the continuous problem in the $H^1$-norm and a mesh-dependent norm, respectively. We present extensive numerical results to illustrate the performance of the approach, including a comparison with the popular Taylor-Hood elements.



\end{abstract}

\section{Introduction}\label{sect:intro}

Topology optimization for fluid mechanics has received a lot of attention in both academia and industry \cite{Aleandersen:2020}. It was initiated by Borrvall and Peterson \cite{BorrvallPetersson:2003}, who propose to minimize the dissipated power plus the potential power of the applied body force (i.e., total potential power) in the Stokes-Darcy flow characterized by a (inverse) variable permeability, so as to find an optimal distribution of fluids corresponding to the Darcy flow in a porous medium and the Stokes flow. Their seminal work has become the foundation of many subsequent investigations; see the references  \cite{Gersborg-Hansen:2005,Olesen:2006,GuestPrevost:2006,WikerKlarbringBorrvall:2007,Age:2008,Kreissl:2011,Alonso:2018,Sa:2018,Deng:2018,Alonso:2020,Thore:2021} for an incomplete list. However, the approach of Borrvall and Petersson only works for the special quadratic energy functional and its extension to a more general objective functional cannot be expected \cite{Evgrafov:2005}. To resolve the issue of the existence of an optimal design, Garcke et al. \cite{GarckeHecht:2016a,GarckeHecht:2016b,GarckeHechtHinzeKahle:2015}  proposed to incorporate a total variation (TV) term into the objective constrained by the stationary (Navier-)Stokes equations, and proposed a phase-field approach using the Ginzburg-Landau energy functional to deal with the numerical challenge.

In this work, we investigate the computational aspect of the phase-field approximation of topology optimization with the Stokes flow. Let $\Omega\subset\mathbb{R}^d$ ($d=2,3$) be an open bounded domain with a polygonal/polyhedral boundary $\partial\Omega$ with a unit outward normal $\bold{n}$. Then the topology optimization problem in the porous-medium phase-field formulation reads \cite{GarckeHecht:2016a,GarckeHecht:2016b,GarckeHechtHinzeKahle:2015}
\begin{subequations}\label{min}
    \begin{align}
        & \qquad \qquad \inf_{\phi \in \U}\mathcal{J}^\eps(\phi) : =  \frac{1}{2}\int_\Omega \alpha_{\eps}(\phi) |\bold{u}|^2 \dx + \int_\Om G(\bold{x},\bold{u},\bold{\nabla}\bold{u}) \dx + \gamma \mathcal{P}_\eps(\phi), \label{min_phase-field}\\
        & \text{subject to}~\bold{u} \in \bold{U}:~	
             \mu \int_\Omega \bold{\nabla} \bold{u} : \bold{\nabla} \bold{v} \dx + \int_\Omega \alpha_\eps(\phi) \bold{u} \cdot \bold{v} \dx = \int_{\Omega}  \bold{f} \cdot \bold{v} \dx, \quad \forall \bold{v} \in \bold{V}. \label{vp_stokes_div-free}
    \end{align}
\end{subequations}
In the objective \eqref{min_phase-field}, $\mathcal{P}_\eps$ (with $\eps>0$ being a relaxation constant) and the admissible set $\mathcal{U}$ of phase-field functions $\phi$ are respectively given by
\begin{equation*}
\mathcal{P}_\eps(\phi) = \frac{\eps}{2}\int_\Om |\bold{\nabla}\phi|^2 \dx + \frac{1}{\eps} \int_\Om f(\phi)  \dx \quad\mbox{and}\quad \U : = \left\{\phi \in H^1(\Omega)~|~\phi \in [0,1]~\text{a.e. in}~\Omega, \int_{\Om}\phi\dx \leq \beta |\Om| \right\} ,
\end{equation*}
where $\beta\in (0,1)$, $\gamma>0$ and $\mu>0$ are the maximum volume fraction that can be occupied by the fluid in the design domain, regularization parameter and fluid viscosity, respectively, $f(\phi)$ is a double obstacle potential $\frac{1}{2}\phi(1-\phi)$ or a double well potential $\frac{1}{4}\phi^2(1-\phi)^2$, and $|\Omega|$ denotes the Lebesgue measure of $\Omega$. The function spaces $\bold{U}$ and $\bold{V}$ are defined by $\bold{U}:=\{\bold{v}\in \bold{H}^1(\Omega)~|~\mathrm{div}\,\bold{v} = 0,\bold{v}|_{\partial\Omega} = \bold{g} \}$ and $\bold{V}:=\{ \bold{v} \in \bold{H}_0^1(\Omega)~|~\mathrm{div}\,\bold{v} = 0\}$, respectively, where the function $\bold{g}$ satisfies the compatibility condition $\int_{\partial\Omega}\bold{g}\cdot\bold{n}\mathrm{d}s = 0$. Note that the inverse permeability tensor $\alpha_\eps$ is coupled with the phase field parameter $\eps$ so that it stays uniformly bounded for any fixed $\eps>0$. The design domain $\Omega$ consists of a fluid region, a non-fluid region occupied by the porous medium and the interfacial thin layer with the width of order $\eps$ represented by a phase-field function $\phi \in \U$ as $\phi \equiv 1$, $\phi \equiv 0$ and $0<\phi<1$ respectively.

Numerically, a gradient flow approach (i.e., the Allen-Cahn or Cahn-Hilliard system) is often employed in the numerical treatment of \eqref{min_phase-field} coupled with the (Navier-)Stokes equations \cite{GarckeHechtHinzeKahle:2015}. Then the resulting problem is approximated by conforming linear elements for the phase field and the lowest-order conforming Taylor-Hood ($P_2$-$P_1$) elements for the velocity and the pressure.  However, even if inf-sup stable finite element (FE) pairs of high order are employed to solve \eqref{vp_stokes_div-free}, the low permeability, corresponding to a large $\alpha$ in the non-fluid region ($\phi \equiv 0$), has an adversary effect on the precision of numerical solutions. Moreover, more degrees of freedom in high-order elements lead to huge computational cost, especially in the three-dimensional case \cite{Thore:2021}. Hence, it is a natural idea to replace the Taylor-Hood elements for the (Navier-)Stokes equations by nonconforming linear (Crouzeix-Raviart or CR) elements \cite{CrouzeixRaviart:1973} and piecewise constants. In this work, we numerically investigate and theoretically analyze the convergence of the CR element for problem \eqref{min_phase-field}-\eqref{vp_stokes_div-free}, and prove the following convergence: Given a sequence of uniformly refined shape-regular meshes $\{\cT_k\}_{k\geq0}$ generated from an initial mesh $\cT_0$, the sequence of discrete minimizing pairs $\{(\phi_k^\ast,\bold{u}^\ast_k)\}_{k\geq 0}$ converges strongly to a minimizing pair to \eqref{min_phase-field}-\eqref{vp_stokes_div-free} up to a subsequence. The key challenge in the analysis is the nonconformity of the CR FE space, for which we establish a novel discrete compactness result.

The overall convergence analysis proceeds as follows. We indicate the dependence on a mesh $\cT_k$ by the subscript $k$. First, we show the existence of a discrete minimizer to problem \eqref{dismin_phase-field}-\eqref{disvp_stokes_div-free} in Theorem \ref{thm:existence_disc}. Then we establish an important weak continuity result of a discrete solution map $\S_k: \phi_k \mapsto \bold{u}_k: = \S_k(\phi_k)$ with respect to the weak $H^1(\Omega)$-topology of $\{\phi_k\}_{k\geq 0}$ in Lemma \ref{lem:vel_weak-conv}. Using the connection operator \cite{Brenner:2003} linking $\bold{X}_k$ to its conforming Lagrange FE space, cf. \eqref{def:con_op}-\eqref{con_operator-stab}, we obtain the $L^2$ strong convergence of $\{\S_k(\phi_k)\}_{k\geq0}$ in Lemma \ref{lem:vel_L2-strong-conv}, in the spirit of discrete compactness associated with nonconforming finite elements in \cite{Stummel:1980} \cite[Chapter 10.3]{ShiWang:2013}. This crucial property allows proving the $L^2$ strong convergence of $\{\bold{\nabla}_k\S_k(\phi_k)\}_{k\geq 0}$ in Lemma \ref{lem:vel_strong-conv}. In view of Lemma \ref{lem:sol-map_cont} and Lemmas \ref{lem:vel_weak-conv}-\ref{lem:vel_strong-conv}, a subsequence of discrete minimizers $\{\phi_{k_j}^\ast\}_{j\geq 0}$ converges strongly in $H^1(\Omega)$ to a minimizer to problem \eqref{min_phase-field}-\eqref{vp_stokes_div-free} and then the associated $\{\bold{u}_{k_j}\}_{j\geq0}$ converges strongly to the corresponding velocity field in Theorem \ref{thm:conv_uniform}. In addition, we show the $L^2$ strong convergence for discrete pressure fields $\{p_{k_j}^\ast\}_{j\geq 0}$ implicit in \eqref{disvp_stokes_div-free} in Theorem \ref{thm:conv_pre_uniform}.

The numerical analysis of topology optimization in fluid flow was studied in \cite{BorrvallPetersson:2003,PapadopoulosSuli:2022,Papadopoulos:2023} for the Borrvall and Petersson model \cite{BorrvallPetersson:2003}. For a piecewise constant FE approximation of the material distribution and an inf-sup stable quadrilateral FE approximation of the velocity and pressure, Borrvall and Petersson \cite{BorrvallPetersson:2003} showed the convergence of the discrete velocity and material distribution to a minimizer weakly in $H^1(\Omega)^d$ and weakly $\ast$ in $L^\infty(\Omega)$, respectively. For any suitable conforming mixed FE space pair such that the velocity and pressure spaces are inf-sup stable, Papadopoulos and Suli \cite{PapadopoulosSuli:2022} proved that for every isolated minimizer, there exists a sequence of FE solutions to the
discretized first-order optimality conditions that strongly converge to the minimizer as the mesh size tends to zero in $H^1(\Omega)^d$ and $L^s(\Omega)$, $s\in [1,\infty)$, respectively. Later Papadopoulos \cite{Papadopoulos:2023} extended the analysis to divergence-free discontinuous Galerkin (DG) methods with an interior penalty, and proved that, given an isolated minimizer of the model, there exists a sequence of DG finite element solutions, satisfying necessary first-order optimality conditions, that strongly converges to the minimizer. In the work, we establish the convergence of CR FE approximations of phase-field dependent topology optimization problem \eqref{min_phase-field}-\eqref{vp_stokes_div-free}, which complements the existing works \cite{BorrvallPetersson:2003,PapadopoulosSuli:2022,Papadopoulos:2023}.

Throughout, we use standard notation for Sobolev spaces \cite{AdamsFournier:2003} and related (semi-)norms. The notation $c$, with or without subscript, denotes a generic positive constant independent of the mesh size and it may take a different value at each occurrence. For a matrix $\bold{\tau}\in\mathbb{R}^{d\times d}$ and a vector $\bold{v}\in\mathbb{R}^d$, $\bold{\tau}\cdot \bold{v}$ is defined row-wisely.
The divergence operator $\mathbf{div}$ is applied to $\bold{\tau}$ row-wisely, i.e.
$\mathbf{div}\,\bold{\tau}=(\sum_{j=1}^d\partial_j\tau_{ij})_{1\leq i\leq d}$. The notation $(\cdot,\cdot)_{L^2(D)}$ denotes the $L^2(D)$ inner product, and we abbreviate it to $(\cdot,\cdot)$ when $D=\Omega$. $L_0^2(\Omega)$ is a subspace of $L^2(\Omega)$ with zero integral average over $\Omega$.

\section{Crouzeix-Raviart FE approximation of topology optimization in fluid}

In this section, we develop the FE approximation of fluid topology optimization.
\subsection{Problem formulation}

First we give the assumptions on problem  \eqref{min_phase-field}-\eqref{vp_stokes_div-free} \cite{GarckeHechtHinzeKahle:2015}.
\begin{assumption}\label{ass:problem}
The following assumptions hold on problem \eqref{min_phase-field}-\eqref{vp_stokes_div-free}.
\begin{itemize}
  \item [(i)] The inverse permeability $\alpha_\eps:[0,1]\to[0,\overline{\alpha}_\eps]$ is decreasing, surjective and continuously differentiable.

  \item [(ii)] The function $G: \Omega \times \R^d \times \R^{d\times d}\to \R$ in the functional $\int_\Omega G(\bold{x},\bold{v}(\bold{x}),\bold{\tau}(\bold{x})) \dx$ is continuous: $\bold{v}_n \to \overline{\bold{v}}$ strongly in $\bold{L}^2(\Omega)$ and $\bold{\tau}_n \to \overline{\bold{\tau}}$ strongly in $\mathbb{L}^2(\Omega):=L^2(\Omega)^{d \times d}$ imply
      \begin{equation}\label{strong-cont}
        \lim_{n\to\infty} ( G(\bold{x},\bold{v}_n,\bold{\tau}_n),1) = ( G(\bold{x},\overline{\bold{v}},\overline{\bold{\tau}}),1).
      \end{equation}
  \item [(iii)] The  functional $\int_\Omega G(\bold{x},\bold{v}(\bold{x}),\bold{\tau}(\bold{x})) \dx$ is weakly lower semi-continuous: $\bold{v}_n \to \overline{\bold{v}}$ strongly in $\bold{L}^2(\Omega)$ and $\bold{\tau}_n \rightarrow \overline{\bold{\tau}}$ weakly in $\mathbb{L}^2(\Omega)$ imply
      \begin{equation}\label{weak-low-semicont}
        \liminf_{n\to\infty}( G(\bold{x},\bold{v}_n,\bold{\tau}_n),1) \geq ( G(\bold{x},\overline{\bold{v}},\overline{\bold{\tau}}),1).
      \end{equation}
      and is bounded from below over $\bold{L}^2(\Omega)\times \mathbb{L}^2(\Omega)$.

  \item [(iv)] The body force $\bold{f}$ in \eqref{vp_stokes_div-free} belongs to $\bold{L}^2(\Omega)$ and the velocity $\bold{g}\in \bold{H}^{1/2}(\partial\Omega)$ satisfies $\int_{\partial\Omega} \bold{g}\cdot\bold{n}\ds = 0$.
\end{itemize}
\end{assumption}

One sufficient condition for (ii) is that $G$ is a Carath\'{e}odory function, i.e. $G(\cdot, \bold{\xi}, \bold{A}): \Omega \to \R $ is measurable for each pair $(\bold{\xi}, \bold{A}) \in \R^d \times \R^{d\times d} $ and $G(\bold{x}, \cdot, \cdot): \R^d \times \R^{d \times d} \to \R$ is continuous for almost every $\bold{x} \in \Omega$, $G$ satisfies a quadratic growth condition: there exist $a\in L^1(\Omega)$ and $b_1,b_2\in L^\infty(\Omega)$ such that for almost every $\bold{x} \in \Omega$,
      \begin{equation}\label{growth_cond}
        \left| G(\bold{x},\bold{\xi},\bold{A}) \right|
        \leq a(\bold{x}) + b_1(\bold{x}) |\bold{\xi}|^{2} + b_2(\bold{x}) |\bold{A}|^2,\quad \forall \bold{\xi} \in \R^d, \forall \bold{A} \in \R^{d \times d}.
      \end{equation}
(ii)-(iii) imply that the functional $\int_\Omega G(\bold{x},\bold{u}(\bold{x}),\bold{\nabla u}(\bold{x})) \dx$ is continuous and weakly lower semi-continuous in $\bold{H}^1(\Omega)$, respectively, and bounded from below over $\bold{U}$. It covers the commonly used $\frac{\mu}{2}\int_{\Omega} |\bold{\nabla}\bold{u}|^2 \dx - \int_\Omega \bold{f} \cdot \bold{u} \dx$ or $\frac{1}{2} \int_\Omega | \bold{u} - \bold{u}_d |^2 \dx$ over $\bold{U}$ with a given $\bold{u}_d \in \bold{U}$. 

In light of (iv), Lax-Milgram lemma guarantees the unique solvability of problem \eqref{vp_stokes_div-free} for each fixed $\phi\in \U$, i.e., a well-defined solution map $\S: \U \to \bold{U}$ \cite[Lemma 1]{GarckeHecht:2016b}.  We have the following continuity of $\S$ with respect to the $L^1(\Omega)$ topology. Then the direct method in calculus of variations \cite{Dacorogna:2008} guarantees the existence of at least one minimizer $\phi^\ast\in\U$ to problem \eqref{min_phase-field}-\eqref{vp_stokes_div-free} \cite[Theorem 1]{GarckeHecht:2016b} \cite[Theorem 2.9]{GarckeHecht:2016a}.
\begin{lemma}\label{lem:sol-map_cont}
Let Assumption \ref{ass:problem} (i) and (iv) hold. Let $\{\phi_n\}_{n\geq 0} \subset \U$ converge to $\phi \in U$ in $L^1(\Omega)$. Then $\{\S(\phi_n)\}_{n\geq 0}$ converges to $\S(\phi)$ strongly in $H^1(\Omega)$.
\end{lemma}
\begin{proof}
By \eqref{vp_stokes_div-free} and direct computation, we obtain
    \[
        \mu (\bold{\nabla}(\S(\phi) - \S(\phi_n)), \bold{\nabla}\bold{v}) +  (\alpha_\eps(\phi_n)(\S(\phi) - \S(\phi_n)), \bold{v}) =   ( (\alpha_\eps(\phi_n) - \alpha_\eps(\phi))\S(\phi) , \bold{v}), \quad \forall  \bold{v} \in \bold{V}.
    \]
    Taking $\bold{v} = \S(\phi) - \S(\phi_n) \in \bold{V}$ and applying Poincar\'{e}'s inequality yield
    \[
        |\S(\phi) - \S(\phi_n)|_{\bold{H}^1(\Omega)} \leq c \|(\alpha_\eps(\phi_n) - \alpha_\eps(\phi))\S(\phi)\|_{\bold{L}^2(\Omega)}.
    \]
The $L^1(\Omega)$ convergence of $\{\phi_n\}_{n\geq 0}$ to $\phi$ implies a.e. convergence up to a subsequence \cite[Theorem 1.21]{EvansGariepy:2015}. By the inequality $|\alpha_\eps(\phi_n) - \alpha_\eps(\phi)|\leq \overline{\alpha}_\eps$, Lebesgue dominated convergence theorem \cite[Theorem 1.19, p. 28]{EvansGariepy:2015} and Poincar\'{e} inequality, we obtain $\|\S(\phi) - \S(\phi_n)\|_{\bold{H}^1(\Omega)}\to 0$ up to a subsequence. By a standard subsequence contradiction argument, the whole sequence converges.
\end{proof}

\subsection{Crouzeix-Raviart FE approximation}
Now we describe the FE approximation of problem \eqref{min_phase-field}-\eqref{vp_stokes_div-free}. Let $\cT$ be a conforming and shape-regular triangulation of the domain $\overline{\Omega}$ into closed triangles/tetrahedra
\cite{Ciarlet:2002}. The collection of all edges/faces (respectively interior edges/faces) in $\cT$ is denoted by $\F_\cT$ (respectively
$\F_{\cT}(\Omega)$). For each $F\in\mathcal{F}_\cT$, we assign a fixed normal unit vector $\bold{n}_F$, which points from
$T_+$ to $T_-$ if two elements $T_+, T_-\in \cT$ share a common $F \in\mathcal{F}_\cT(\Omega)$, and is chosen to be the
unit outward normal $\bold{n}$ to $\Omega$ if $F\subset\partial\Omega$, and $\bold{n}_{\partial T}$ denotes the unit outward normal
to the boundary of each $T\in\cT$. $\mathcal{N}_{\cT}(F)$ denotes the set of elements sharing an $F\in \mathcal{F}_\cT$
and $\mathcal{N}_{\cT}(T)$ the set of an $T\in \cT$ and its neighboring elements in $\cT$, the union of which is denoted
by $\omega_\cT(T)$, respectively. The local mesh-size $h_T$ (respectively $h_F$) of each $T\in \cT$ (respectively $F\in \F_\cT$) is defined
as $|T|^{1/d}$ (respectively $|F|^{1/(d-1)}$), which further yields a piecewise constant mesh-size function $h_\cT:\overline{\Omega}
\to \mathbb{R}^+$ with $h_\cT|_{T} : = |T|^{1/d}$ for each $T\in \cT$. $P_m(T)$ denotes the set of all polynomials of total degree no more than $m$ on $T \in \cT$ for any given integer $m\geq0$. Now over $\cT$, we define the discrete admissible set $\U_\cT:=S_\cT\cap\U$ ($S_\cT$ is the usual $H^1(\Omega)$-conforming linear FE space), the nonconforming linear FE (CR) space \cite{CrouzeixRaviart:1973} and the piecewise constant space, respectively, by
\begin{align}\label{C-R_space}
    \bold{X}_\cT &: = \left\{ \bold{v}_\cT \in \bold{L}^2(\Omega)~|~\bold{v}_\cT|_T \in P_1(T)^d~\forall T\in\cT,~\int_{F}[\bold{v}_\cT]\ds=\bold{0}~\forall F\in\mathcal{F}_{\cT}(\Omega)\right\},\\
\label{discretespace_pressure}
    Q_{\cT}&:=\{v_\cT \in L_{0}^{2}(\Omega)~|~v_{\cT}|_T \in P_{0}(T)~\forall T\in\mathcal{T}\},
\end{align}
where $[\boldsymbol{v}_\cT]:= \bold{v}_\cT|_{F\subset T_+} - \bold{v}_\cT|_{F\subset T_-}$ is the jump across an interior face/edge $F = \partial T_+ \cap \partial T_-$ and takes the value of
$\boldsymbol{v}_\cT$ on $F\subset\partial\Omega$. Any function in $\bold{X}_\cT$ is determined by its values at the centers of faces/edges in $\F_\cT$ and those members in $\bold{X}_\cT$ vanishing at the boundary centers comprise the subspace $\bold{Z}_\cT$. With $\bold{\nabla}_\cT$ (respectively $\mathrm{div}_\cT$) denoting the piecewise gradient (respectively divergence) operator over $\cT$, we have $\mathrm{div}_\cT \bold{v}_\cT \in Q_\cT$ for any $\bold{v}_\cT \in \bold{Z}_\cT$. Thus one can define a piecewise divergence free finite element subspace
\begin{equation}\label{C-R_space_div-free}
    \bold{V}_\cT : = \left\{\bold{v}_\cT \in \bold{Z}_\cT ~|~(\mathrm{div}_{\cT}  \bold{v}_\cT , q_\cT) = 0~\forall q_\cT \in Q_\cT\} = \{\bold{v}_\cT \in \bold{Z}_\cT ~|~\mathrm{div}_{\cT} \bold{v}_\cT = 0\right\}.
\end{equation}
Also on the CR FE space $\bold{X}_\cT$, we define an interpolation operator $\bold{\Pi}_\cT:\bold{H}^1(\Omega) \to \bold{X}_\cT$ such that \cite{CrouzeixRaviart:1973, OrtnerPraetorius:2011}
\begin{equation}\label{int_operator-cr}
  \int_F \bold{\Pi}_\cT \bold{v} \dd s = \int_F \bold{v} \dd s, \quad \forall F \in \mathcal{F}_\cT.
\end{equation}
It satisfies the first order approximation property
\begin{equation}\label{int_operator-cr_approx}
  \|\bold{v}-\bold{\Pi}_\cT\bold{v}\|_{\bold{L}^2(T)}\leq c h_T \|\bold{\nabla}\bold{v}\|_{\bold{L}^2(T)}, \quad \forall T \in \cT,
\end{equation}
with positive constant $c$ depending on the shape regularity of $\cT$ and the stability estimate
\begin{equation}\label{int_operator-cr_stab}
    \|\bold{\nabla}(\bold{\Pi}_\cT \bold{v})\|_{\bold{L}^2(T)} \leq \|\bold{\nabla}\bold{v}\|_{\bold{L}^2(T)}.
\end{equation}
By Assumption \ref{ass:problem}(iv) and \cite[Lemma 2.2, Chapter I]{GiraultRaviart:1986}, the boundary condition $\bold{g}$ admits a lift $\bold{w}\in \bold{U}$.
{Moreover}, by \eqref{int_operator-cr} and $\mathrm{div}\,\bold{w}=0$, {there holds} $\mathrm{div}_\cT \bold{\Pi}_\cT\bold{w} = 0$.
The FE approximation for \eqref{min_phase-field}-\eqref{vp_stokes_div-free} reads
\begin{subequations}\label{dismin}
    \begin{align}
         &\qquad\qquad\inf_{\phi_\cT \in \U_\cT}\mathcal{J}^\eps_\cT(\phi_\cT) : =  \frac{1}{2}\int_\Om \alpha_{\eps}(\phi_\cT)
         |\bold{u}_\cT|^2 \dx + \int_\Om G(\bold{x},\bold{u}_\cT,\bold{\nabla}_\cT\bold{u}_\cT) \dx + \gamma \mathcal{P}_\eps(\phi_\cT),
         \label{dismin_phase-field}\\
        & \text{s.t.}~\bold{u}_\cT \in \bold{\Pi}_\cT \bold{w} + \bold{V}_\cT:~	
             \mu (\bold{\nabla}_\cT \bold{u}_\cT, \bold{\nabla}_\cT \bold{v}_\cT) + ( \alpha_\eps(\phi_\cT) \bold{u}_\cT, \bold{v}_\cT)  =  (\bold{f}, \bold{v}_\cT), \quad \forall  \bold{v}_\cT \in \bold{V}_\cT. \label{disvp_stokes_div-free}
    \end{align}
\end{subequations}
The existence of a minimizer $\phi_\cT^*$ to the discrete problem is given in Theorem \ref{thm:existence_disc}.

\subsection{Numerical solver}\label{subsect:solver}
By the augmented Lagrangian method for the volume constraint, the phase-field model for \eqref{min_phase-field}-\eqref{vp_stokes_div-free} is given by
\begin{equation}\label{eq:Lagrangian}
\begin{aligned}
\min _{\phi \in[0,1]} \mathcal{L}(\phi,\bold{u}) = \mathcal{J}^\eps(\phi)+\ell W(\phi)+\frac{\zeta}{2}W(\phi)^2\quad
\text {subject to}  \quad (\ref{min}\text{b}),
\end{aligned}
\end{equation}
where $W(\phi):=\int_\Omega \phi \dx-\beta|\Omega|$, $\ell$ denotes a Lagrange multiplier, $\beta\in(0,1)$ represents a volume fraction, and $\zeta>0$ is a penalty parameter. We adopt the gradient flow approach to solve \eqref{eq:Lagrangian} \cite{TakezawaNishiwakiKitamura:2010,QHZ:2022,JinLiXuZhu:2024,LiXuZhu:2023,HQZ:2023}. By introducing the pseudo-time $t>0$, the phase-field equation is given by
\begin{equation}\label{eq:phasefield_evolve}
\left\{\begin{aligned}
\partial_t \phi &=-\frac{\delta \mathcal{L}(\phi, \bold{u})}{\delta \phi}, && \mbox{in }\Omega, \, t>0, \\
\phi(\cdot, 0)&=\phi^0, && \mbox{in }\Omega, \\
\partial_{\bold{n}} \phi&=0, && \mbox{on }\partial \Omega, t>0,
\end{aligned}\right.
\end{equation}
where $\phi^0\in L^2(\Omega)$ is the initial guess. The function $f(\phi)$ in $\mathcal{J}^\eps(\phi)$ is chosen to be a double-well potential $\frac{1}{4}\phi^2(1-\phi)^2$ and $\alpha_{\varepsilon}(\phi)$ is taken to be $\alpha_{0}(1-\phi)^2$ with $\alpha_{0}>0$.

{We employ a nested inner-outer iterative algorithm (see Algorithm \ref{Alg:UniformRefinement}) for \eqref{eq:phasefield_evolve} with $m=1,2,\cdots,M$ and $n=1,2,\cdots,N$ denoting the inner and outer iteration indices, respectively. The discretization of \eqref{eq:phasefield_evolve} in the temporal direction}  
is to find $\phi^{n, m+1}\in H^1(\Omega)$ such that \cite{LiYang:2022}
\begin{align}
& (\Delta t)^{-1}(\phi^{n, m+1}, \psi)+\eps \gamma(\nabla \phi^{n, m+1}, \nabla \psi) +((\tfrac{1}{2}\alpha_0|\bold{u}^{n}|^2+\tilde{S}) \phi^{n, m+1}, \psi) \label{eq:phasefieldupdateG1}\\
= & (\Delta t)^{-1}(\phi^{n, m}, \psi)+(-\eps^{-1}\gamma f(\phi^{n, m})+\alpha_0|\bold{u}^{n}|^2-\ell^n-\zeta^{n} W(\phi^{n,m}), \psi)  ((\tilde{S}-\tfrac{1}{2}\alpha_0|\bold{u}^{n}|^2) \phi^{n, m}, \psi),
\quad \forall\psi\in H^1(\Omega),\nonumber
\end{align}
where $\tilde{S}$ is a stabilized parameter, and $\Delta t$ is the time step size.
The Lagrange multiplier $\ell$ is updated using a Uzawa type scheme by
\begin{equation}\label{eq:multiplier_update}
\ell^{n+1} = \ell^{n}+{\zeta}^{n} W(\phi^{n+1}).
\end{equation}
The value of $\ell$ increases (respectively decreases) as the current structure volume is larger (respectively smaller) than the target volume $\beta |\Omega|$. The penalty
parameter is updated during iterations via
\begin{equation}\label{eq:penalty_update}
\zeta^{n+1}=\kappa\zeta^{n},
\end{equation}
with $\kappa \geq 1$. Numerically, the function $\phi^{n,m}$ may do not preserve the box constraint $[0,1]$. To remedy the issue, we perform a projection step
\begin{equation}\label{phase_field_projection}
\phi^{n,m}\leftarrow\min\{\max\{0,\phi^{n,m}\},1\}.
\end{equation}
Then we set $\phi^{n+1,m}:=\phi^{n,M}$ after $M$ steps of \eqref{eq:phasefield_evolve}.

\begin{algorithm}
    \caption{Phase-field topology optimization}\label{Alg:UniformRefinement}
    \LinesNumbered
    \KwIn{given $\ell^{0}_{0}$, $\zeta^{0}_{0}$, $\kappa$, $\tilde S$, $\Delta t$, $\alpha_{0}$, initial guess $\phi_{0}^{0,0}$, and integers $K$, $N$ and $M$}
    \KwOut{the phase--field function $\phi^{N,M}_{K}$}
    \For {$k=0:K$}{
    \For {$n=1:N$}{
        Solve the state problem \eqref{min_phase-field}-\eqref{vp_stokes_div-free} to obtain $\bold{u}^{n}_{k}$; \\
        \For {$m=1:M$}{
                Update $\phi^{n,m}_{k}$ using \eqref{eq:phasefieldupdateG1}; \\
                Projection $\phi^{n,m}_{k}$ to $[0,1]$ via \eqref{phase_field_projection};
        }
        Update the Lagrange multiplier $\ell^{n}_{k}$ by \eqref{eq:multiplier_update} and penalty parameter $\zeta^{n}_{k}$ by \eqref{eq:penalty_update};
    }
    Refine the mesh uniformly;
    }
\end{algorithm}

\section{Convergence analysis of the FE approximation}\label{sect:conv_uniform}
First we show the existence of a minimizer $\phi_\cT^*\in\mathcal{U}_\cT$ to problem \eqref{dismin_phase-field}-\eqref{disvp_stokes_div-free} on the mesh $\cT$. We equip the FE space $\bold{X}_{\cT}$ with a mesh-dependent semi-norm
$\|\boldsymbol{v}\|^{2}_{1,\cT}:=\sum_{T\in\cT}\int_{T}|\bold{\nabla}\bold{v}|^2\dx=\int_\Omega |\boldsymbol{\nabla}_{\cT}\boldsymbol{v}|^2 \dx$, which is a norm over $\bold{Z}_\cT$ and coincides with the usual norm over $\bold{H}^1_0(\Omega)$ due to Poincar\'{e} inequality. Also we have the following discrete Poincar\'{e} inequality \cite{Brenner:2003, BuffaOrtner:2009}
\begin{equation}\label{disc_Poin_ineq}
    \|\bold{v}_\cT\|_{\bold{L}^2(\Omega)} \leq c_{\text{dp}} \|\bold{v}_\cT\|_{1,\cT},\quad\forall \bold{v}_\cT \in \bold{Z}_\cT,
\end{equation}
with the constant $c_{\text{dp}}$ depending only on the shape regularity of the mesh $\cT$. By \eqref{disc_Poin_ineq} and Lax-Milgram lemma, there exists a unique $\widetilde{\bold{u}}_\cT \in \bold{V}_\cT$ such that
\begin{equation}\label{disvp-homogeneous_stokes_div-free}
    \begin{aligned}
    \mu (\bold{\nabla}_\cT \widetilde{\bold{u}}_\cT, \bold{\nabla}_\cT \bold{v}_\cT) & +  (\alpha_\eps(\phi_\cT) \widetilde{\bold{u}}_\cT, \bold{v}_\cT)  = (\bold{f}, \bold{v}_\cT) \\
    & - \mu (\bold{\nabla}_\cT \bold{\Pi}_\cT \bold{w}, \bold{\nabla}_\cT \bold{v}_\cT) - (\alpha_\eps(\phi_\cT) \bold{\Pi}_\cT \bold{w}, \bold{v}_\cT), \quad \forall  \bold{v}_\cT \in \bold{V}_\cT.
    \end{aligned}
\end{equation}
Further, taking $\bold{v}_\cT = \widetilde{\bold{u}}_\cT$ in \eqref{disvp-homogeneous_stokes_div-free} and noting the stability of $\bold{\Pi}_\cT$ from \eqref{int_operator-cr_approx}-\eqref{int_operator-cr_stab} and Assumption \ref{ass:problem}(i) on $\alpha_\eps$, we have
\begin{equation}\label{stab_discvp-hom}
    \|\widetilde{\bold{u}}_\cT\|_{1,\cT} \leq c \left (\|\bold{f}\|_{\bold{L}^2(\Omega)} + \|\bold{w}\|_{\bold{H}^1(\Omega)}\right).
\end{equation}
The unique solvability of problem \eqref{disvp-homogeneous_stokes_div-free} implies that problem \eqref{disvp_stokes_div-free} has a unique $\bold{u}_\cT = \widetilde{\bold{u}}_\cT + \bold{\Pi}_\cT \bold{w} \in \bold{\Pi}_\cT \bold{w} + \bold{V}_\cT$ for any $\phi_\cT \in \U_\cT$. This also allows defining a discrete solution map $\bold{\mathcal{S}}_\cT: \mathcal{U}_\cT \to \bold{\Pi}_\cT \bold{w} + \bold{V}_\cT$ as $\phi_\cT \mapsto \bold{u}_\cT:= \bold{\mathcal{S}}_\cT(\phi_\cT)$ for problem \eqref{disvp_stokes_div-free}. Further, it follows from \eqref{int_operator-cr_approx}, \eqref{int_operator-cr_stab}, \eqref{disc_Poin_ineq} and \eqref{stab_discvp-hom} that
\begin{equation}\label{stab_discvp}
    \|\S_\cT(\phi_\cT)\|_{\bold{L}^2(\Omega)} +  \|\bold{\nabla}_\cT\S_\cT(\phi_\cT)\|_{\bold{L}^2(\Omega)} \leq c \left (\|\bold{f}\|_{\bold{L}^2(\Omega)} + \|\bold{w}\|_{\bold{H}^1(\Omega)}\right),
\end{equation}
with $c$ depending only on the shape regularity of the mesh $\cT$.

\begin{theorem}\label{thm:existence_disc}
Let Assumption \ref{ass:problem}  hold. Then for any $\eps>0$, there exists one minimizer $\phi^*_\cT$ to \eqref{dismin_phase-field}-\eqref{disvp_stokes_div-free}.
\end{theorem}

\begin{proof}
The proof is similar to that of  \cite[Theorem 2.4]{GarckeHecht:2016a} and \cite[Theorem 1]{GarckeHecht:2016b}. By Assumption \ref{ass:problem}(iii), $\J_\cT^\eps(\phi_\cT)$ is bounded from below. Thus there exist a minimizing sequence $\{\phi^n_\cT\}_{n\geq 0}\subset \U_\cT$ and a sequence of solutions $\{ \S_\cT(\phi^n_\cT)\}_{n\geq0}\subset \bold{\Pi}_\cT \bold{w} + \bold{V}_\cT$ such that $\lim_{n\to\infty}\mathcal{J}_\cT^\eps(\phi_\cT^n) = \inf_{\phi_\cT \in \U_\cT}\mathcal{J}^\eps_\cT(\phi_\cT).$ This and the box constraint $0\leq \phi_\cT^n \leq 1$ in $\Omega$ (cf. $\U_\cT$) imply that $\left\{\|\phi_\cT^n\|_{H^1(\Omega)}\right\}_{n\geq0}$ is uniformly bounded. Since $\widetilde{\bold{u}}_\cT^n: = \S_\cT(\phi_\cT^n) - \bold{\Pi}_\cT \bold{w} \in \bold{V}_\cT$ solves \eqref{disvp-homogeneous_stokes_div-free} with $\phi_\cT$ replaced by $\phi_\cT^n$, it follows from \eqref{disc_Poin_ineq} and \eqref{stab_discvp-hom} that $ \| \widetilde{\bold{u}}_\cT^n \|_{\bold{L}^2(\Omega)} + \|\widetilde{\bold{u}}_\cT^n\|_{1,\cT} \leq c$ for $n\geq 0$. Since $\U_\cT$ is closed in a finite-dimensional space $S_\cT$ and $\bold{V}_\cT$ is finite-dimensional, by the norm equivalence in a finite-dimensional space, we may extract two subsequences, still denoted by $\{\phi_\cT^n , \widetilde{\bold{u}}_\cT^n\}_{n\geq0}$, and find a pair $(\phi^\ast_\cT,\widetilde{\bold{u}}^\ast_\cT) \in \U_\cT \times \bold{V}_\cT$ such that
$$\lim_{n\to\infty}\|\phi_\cT^n - \phi_{\cT}^\ast\|_{H^1(\Omega)} = 0\quad \mbox{and}\quad \lim_{n\to\infty}\|\widetilde{\bold{u}}_\cT^n - \widetilde{\bold{u}}_\cT^\ast \|_{\bold{L}^2(\Omega)} + \|\widetilde{\bold{u}}_\cT^n - \widetilde{\bold{u}}_\cT^\ast \|_{1,\cT} = 0.$$ Thus there exist two subsequences $\left\{\phi_\cT^n\right\}_{n\geq0}$ and $\left\{ \S_\cT(\phi_\cT^n) \right\}_{n\geq0}\subset \bold{\Pi}_\cT \bold{w} + \bold{V}_\cT$ such that
    \begin{align}\label{thm:existence_disc_pf01}
        \phi_\cT^n \to \phi_{\cT}^\ast \quad \text{a.e. in}~\Omega,\\
\label{thm:existence_disc_pf02}
      \S_\cT(\phi_\cT^n) \to  \widetilde{\bold{u}}^\ast_\cT + \bold{\Pi}_\cT\bold{w} \quad \text{in}~ \bold{L}^2(\Omega),\quad \bold{\nabla}_\cT \S_\cT(\phi_\cT^n)\to \bold{\nabla}_\cT (\widetilde{ \bold{u}}^\ast_\cT + \bold{\Pi}_\cT\bold{w} )\quad \text{ in}~\mathbb{L}^2(\Omega).
    \end{align}
    By the relation \eqref{thm:existence_disc_pf01}, Assumption \ref{ass:problem}(i) on $\alpha_\eps$ and Lebesgue dominated convergence theorem, we obtain
    \begin{equation}\label{eqn:conv0}
         \lim_{n\to\infty}\|\alpha_\eps (\phi_\cT^n) \bold{v}_\cT- \alpha_\eps (\phi_\cT^\ast) \bold{v}_\cT\|_{ \bold{L}^2(\Omega)}=0,\quad \forall \bold{v}_\cT\in \bold{V}_\cT.
    \end{equation}
By letting $\bold{u}_\cT^\ast: = \widetilde{\bold{u}}_\cT^\ast + \bold{\Pi}_\cT \bold{w}$ and passing to the limit $n\to\infty$ in \eqref{disvp_stokes_div-free}, it follows from \eqref{thm:existence_disc_pf02} and \eqref{eqn:conv0} that
    \[
        \mu (\bold{\nabla}_\cT \bold{u}_\cT^\ast, \bold{\nabla}_\cT \bold{v}_\cT) + ( \alpha_\eps(\phi_\cT^\ast) \bold{u}^\ast_\cT, \bold{v}_\cT)  = (\bold{f}, \bold{v}_\cT), \quad \forall  \bold{v}_\cT \in \bold{V}_\cT.
    \]
That is, $\bold{u}^\ast_\cT = \S_\cT(\phi_\cT^\ast)$. Next we prove that $\phi_\cT^\ast$ is a minimizer to problem \eqref{dismin_phase-field}-\eqref{disvp_stokes_div-free}. We use the splitting
    \[
       (\alpha_\eps (\phi_\cT^n), |\S_\cT(\phi_\cT^n) |^2 ) =
        (\alpha_\eps (\phi_\cT^n),  |\S_\cT(\phi_\cT^n) |^2 - |\S_\cT(\phi_\cT^\ast) |^2) + ( \alpha_\eps (\phi_\cT^n), |\S_\cT(\phi_\cT^\ast) |^2).
    \]
    Assumption \ref{ass:problem}(i) on $\alpha_\eps$ and the $\bold{L}^2(\Omega)$ convergence in \eqref{thm:existence_disc_pf02} imply
    \[
       \lim_{n\to\infty} \left|(\alpha_\eps (\phi_\cT^n), \left( |\S_\cT(\phi_\cT^n) |^2 - |\S_\cT(\phi_\cT^\ast) |^2 \right))\right|  \leq \lim_{n\to\infty} c \| \S_\cT(\phi_\cT^n) - \S_\cT(\phi_\cT^\ast) \|_{\bold{L}^2(\Omega)} = 0,
    \]
and Lebesgue dominated convergence theorem yields
$       \lim_{n\to\infty} (\alpha_\eps (\phi_\cT^n),|\S_\cT(\phi_\cT^\ast) |^2) = (\alpha_\eps (\phi_\cT^\ast), |\S_\cT(\phi_\cT^\ast) |^2).$
    Therefore,
$\lim_{n\to\infty}   (\alpha_\eps (\phi_\cT^n), |\S_\cT(\phi_\cT^n) |^2 )= (\alpha_\eps (\phi_\cT^\ast), |\S_\cT(\phi_\cT^\ast) |^2)$.
{By Assumption \ref{ass:problem}(ii) and \eqref{thm:existence_disc_pf02}, there holds
    \[
        (G(\bold{x}, \S_\cT(\phi_\cT^\ast), \bold{\nabla}_\cT \S_\cT(\phi_\cT^\ast)),1)  = \lim_{n\to\infty} (G(\bold{x}, \S_\cT(\phi_\cT^n), \bold{\nabla}_\cT \S_\cT(\phi_\cT^n)),1).
    \]}
    The $H^1(\Omega)$ convergence of $\left\{\phi_\cT^n\right\}_{n\geq0}$, \eqref{thm:existence_disc_pf01} and Lebesgue dominated convergence theorem imply
    \[
       \lim_{n\to\infty} \tfrac{\eps}{2}\|\bold{\nabla}\phi^n_\cT\|_{L^2(\Omega)}^2 + \tfrac{1}{\eps} (f(\phi^n_\cT),1) =
        \tfrac{\eps}{2}\|\bold{\nabla}\phi^\ast_\cT\|_{L^2(\Omega)}^2 + \tfrac{1}{\eps}(f(\phi^\ast_\cT),1).
    \]
Noting that $\{\phi_\cT^n\}_{n\geq 0}$ is a minimizing sequence and collecting the above three results lead to the desired assertion.
\end{proof}

Next, we analyze the convergence. The dependence on a triangulation $\cT_k$ is indicated by the refinement level $k\geq 0$ in the subscript.
The main result is the $H^1(\Omega)$ strong convergence of discrete minimizers $\{\phi_k^\ast\}_{k\geq 0} \subset \bigcup_{k\geq 0} \U_k$ (up to a subsequence), over a sequence of conforming, uniformly refined, and shape-regular meshes $\{\cT_k\}_{k\geq0}$ to a minimizer $\phi^*\in \mathcal{U}$ to
problem \eqref{min_phase-field}-\eqref{vp_stokes_div-free}. We employ the direct method in the calculus of variations as in Theorem \ref{thm:existence_disc}. However, it requires the sequential continuity of $\S_k$ along a sequence of discrete phase
fields $\{\phi_k\in\U_k\}_{k\geq0}$ converging to some $\phi \in \U$ instead of the discrete set $\U_k$ with respect to the
weak $H^1(\Omega)$-topology. The nonconformity of the CR FE space, i.e. $\bold{X}_k \nsubseteq \bold{H}^1(\Omega)$, is the major
obstacle. We resort to the following argument:
\begin{itemize}
    \item
The stability of $\bold{u}_k = \S_k(\phi_k)$,
elementwise integration by parts and the vanishing mesh size allow extracting a weakly convergent subsequence, still denoted by
$\{\bold{u}_k\}_{k=1}^\infty$ and some $\overline{\bold{u}} = \S (\phi) \in \bold{U}$ such that $\bold{u}_k \rightarrow \overline{\bold{u}}$ weakly in
$\bold{L}^2(\Omega)$ and $\bold{\nabla}_k\bold{u}_k \to \bold{\nabla}\overline{\bold{u}}$ weakly in $\mathbb{L}^{2}(\Omega)$;
see Lemma \ref{lem:vel_weak-conv}.
\item We upgrade the $\bold{L}^2(\Omega)$ weak convergence of $\{\bold{u}_k\}_{k\geq0}$ to a strong one in Lemma \ref{lem:vel_L2-strong-conv}. This discrete compactness property extends Rellich compactness theorem to
sequences of piecewise polynomials in nonconforming and hybrid FE spaces, and was first discovered by Stummel \cite{Stummel:1980}. The main tool in  \cite{Stummel:1980} is the $L^p$ ($1\leq p < \infty$) equi-continuity of bounded sequences,
also known as Riesz-Tamarkin theorem or Kolmogorov-Riesz-Fr\'{e}chet theorem  \cite[Theorem 2.32, p. 38]{AdamsFournier:2003}.
The proof of Lemma \ref{lem:vel_L2-strong-conv} uses a connection/enriched operator \cite{Brenner:2003} from $\bold{X}_k$ to
its conforming counterpart: vectorial $P_2$ ($d=2$) or $P_3$ ($d=3$) Lagrange FE space.
\item Under the assumption
$\|\bold{\nabla}_k(\bold{w} - \bold{\Pi}_k\bold{w})\|_{\bold{L}^2(\Omega)}\to 0$, the compactness enables deriving a strong convergence
$\|\bold{\nabla}_k(\bold{u}_k - \overline{\bold{u}})\|_{\bold{L}^2(\Omega)}\to 0$ in Lemma \ref{lem:vel_strong-conv}.
\end{itemize}

\begin{lemma}\label{lem:vel_weak-conv}
Let Assumption \ref{ass:problem}(i) and (iv) hold. Let $\{\phi_k\in\U_k\}_{k\geq 0}$ converge to some $\phi \in \U $ weakly in $H^1(\Omega)$ as $\lim_{k\to\infty}\|h_k\|_{L^\infty(\Omega)}= 0$. Then the sequence of discrete velocities $\{\S_k(\phi_k)\in \bold{\Pi}_k \bold{w} + \bold{V}_k\}_{k\geq0}$ converges to the unique solution  $\S(\phi) \in \bold{U}$ of problem \eqref{vp_stokes_div-free}:
\begin{equation}\label{vel_weak-conv}
        \S_k (\phi_k) \rightarrow \S (\phi)\quad \text{weakly in}~\bold{L}^2(\Omega), \quad \bold{\nabla}_k \S_k (\phi_k) \rightarrow \bold{\nabla} \S (\phi) \quad \text{weakly in}~\mathbb{L}^2(\Omega).
    \end{equation}
\end{lemma}

\begin{proof}
Since $\widetilde{\bold{u}}_k \in \bold{V}_k$ solves problem \eqref{disvp-homogeneous_stokes_div-free} in $\cT_k$ and $\bold{\Pi}_k: \bold{H}^1(\Omega) \to \bold{X}_k$, it follows from \eqref{int_operator-cr_approx}, \eqref{int_operator-cr_stab}, \eqref{disc_Poin_ineq} and \eqref{stab_discvp-hom} that both $\{\|\widetilde{\bold{u}}_k\|_{\bold{L}^2(\Omega)} + \|\widetilde{\bold{u}}_k\|_{1,k}\}_{k\geq 0}$ and $\{\|\bold{\Pi}_k \bold{w}\|_{\bold{L}^2(\Omega)} + \|\bold{\Pi}_k \bold{w}\|_{1,k}\}_{k\geq0}$ are uniformly bounded. This implies that
$\{ \|\S_k(\phi_k)\|_{\bold{L}^2(\Omega)} + \|\bold{\nabla}_k\S_k(\phi_k)\|_{\bold{L}^2(\Omega)}\}$ is uniformly bounded.
Thus there exist two subsequences, denoted by $\{\phi_k\}_{k\geq 0}$ and $\{\S_k(\phi_k)\}_{k\geq 0}$, and some $\overline{\bold{u}}\in \bold{L}^2(\Omega)$  and $\bold{\sigma} \in \mathbb{L}^2(\Omega)$ such that
    \begin{equation}\label{lem:vel_weak-conv_pf01}
        \phi_{k} \to \phi \quad \text{a.e. in}~\Omega,\quad
        \S_k(\phi_k) \rightarrow \overline{\bold{u}}  \quad \text{weakly in}~\bold{L}^2(\Omega),\quad  \bold{\nabla}_k\S_k(\phi_k) \rightarrow \bold{\sigma} \quad \text{weakly in}~\mathbb{L}^2(\Omega)\quad \text{as}~k\to\infty.
    \end{equation}
The remaining proof is devoted to the following four assertions: (1) $\overline{\bold{u}}\in \bold{H}^1(\Omega)$ with $\bold{\sigma} = \bold{\nabla}\overline{\bold{u}}$; (2) $\overline{\bold{u}} = \bold{g}$ on $\partial\Omega$; (3) $\mathrm{div}\, \overline{\bold{u}} = 0$; (4) $\overline{\bold{u}}=\S(\phi)$. We divide the lengthy proof into four steps.

    \noindent \textit{Step 1.} $\overline{\bold{u}}\in \bold{H}^1(\Omega)$ with $\bold{\sigma} = \bold{\nabla}\overline{\bold{u}}$. By elementwise integration by parts, we obtain for any $\boldsymbol{\varphi} \in C^\infty(\overline{\Omega})^{d\times d}$,
    \begin{align}\label{lem:vel_weak-conv_pf02}
        ( \S_k(\phi_k), \mathbf{div}\, \bold{\varphi})& =
         \sum_{T\in\cT_{k}}\bigg( - (\bold{\nabla} \S_k(\phi_k) ,\bold{\varphi})_{L^2(T)} + \sum_{F\subset \partial T} (\S_k(\phi_k), \bold{\varphi} \cdot \bold{n}_{\partial T})_{L^2(F)} \bigg)\nonumber\\
        & = - \sum_{T\in\cT_{k}}(\bold{\nabla} \S_k(\phi_k) , \bold{\varphi})_{L^2(T)} + \sum_{F\in\mathcal{F}_{k}}([\S_k(\phi_k)], \bold{\varphi}\cdot\bold{n}_F)_{L^2(F)}.
    \end{align}
 Since $\S_k(\phi_k)\in \bold{\Pi}_k \bold{w} + \bold{V}_k$, we have
 \begin{align*}
     \int_{F} [\S_k(\phi_k)] \dd s &= \bold{0}, \quad \forall F\in \mathcal{F}_k(\Omega),\\ \int_F \S_k(\phi_k) \dd s &= \int_{F} \bold{\Pi}_k \bold{w} \dd s = \int_F \bold{w}  \dd s,\quad \forall F\in \mathcal {F}_k \setminus \mathcal{F}_k(\Omega).
\end{align*}
Thus inserting $\bold{\varphi}_F = \int_F \bold{\varphi} \dd s / | F |$, $ [\S_k(\phi_k)]_F = \int_F [\S_k(\phi_k)] \dd s / |F|$ for $F\in \mathcal{F}_k(\Omega) $ and $[\S_k(\phi_k) - \bold{\Pi}_k\bold{w}]_F = \int_F (\S_k(\phi_k) - \bold{\Pi}_k\bold{w}) \dd s / |F| $ for $F \in \mathcal{F}_k \setminus \mathcal{F}_k(\Omega)$ in the second summand in \eqref{lem:vel_weak-conv_pf02} gives
        \begin{align}\label{lem:vel_weak-conv_pf03}
        & \quad \sum_{F\in\mathcal{F}_{k}} ([\S_k(\phi_k)], (\bold{\varphi}\cdot\bold{n}_F)) _{L^2(F)} \nonumber \\
        &= \sum_{{F\in\mathcal{F}_{k}}} ([\S_k(\phi_k)], (\bold{\varphi}-\bold{\varphi}_F)\cdot\bold{n}_F)_{L^2(F)}  + \sum_{F\in\mathcal{F}_k\setminus\mathcal{F}_{k}(\Omega)} (\bold{\Pi}_k \bold{w}, \bold{\varphi}_F\cdot\bold{n})_{L^2(F)} \nonumber  \\
        & = \sum_{F\in\mathcal{F}_{k}(\Omega)}( [\S_k(\phi_k)], (\bold{\varphi}-\bold{\varphi}_F)\cdot\bold{n}_F)_{L^2(F)} + ( \bold{w} , \bold{\varphi} \cdot \bold{n} )_{L^2(\partial\Omega)}  \nonumber  \\
        & \quad  + \sum_{F\in\mathcal{F}_k\setminus\mathcal{F}_{k}(\Omega)}(\S_k(\phi_k) - \bold{\Pi}_k \bold{w}, (\bold{\varphi}-\bold{\varphi}_F)\cdot\bold{n})_{L^2(F)}  + \sum_{F\in\mathcal{F}_k\setminus\mathcal{F}_{k}(\Omega)} ( \bold{\Pi}_k \bold{w} - \bold{w}, (\bold{\varphi}-\bold{\varphi}_F)\cdot\bold{n})_{L^2(F)} \nonumber \\
        & = \sum_{F\in\mathcal{F}_{k}(\Omega)} \underbrace{([\S_k(\phi_k)] - [\S_k(\phi_k)]_F , (\bold{\varphi}-\bold{\varphi}_F)\cdot\bold{n}_F)_{L^2(F)}}_{{\rm I}_F} + \sum_{F\in\mathcal{F}_k\setminus\mathcal{F}_{k}(\Omega)} \underbrace{(\bold{\Pi}_k \bold{w} - \bold{w}, (\bold{\varphi}-\bold{\varphi}_F)\cdot\bold{n}) _{L^2(F)}}_{{\rm II}_F}  \nonumber \\
        & \quad + \sum_{F\in\mathcal{F}_k\setminus\mathcal{F}_{k}(\Omega)} \underbrace{(\S_k(\phi_k) - \bold{\Pi}_k \bold{w} - [\S_k(\phi_k) - \bold{\Pi}_k\bold{w}]_F , (\bold{\varphi}-\bold{\varphi}_F)\cdot\bold{n} )_{L^2(F)}}_{{\rm III}_F}+ (\bold{w}, \bold{\varphi} \cdot \bold{n})_{L^2(\partial\Omega)}.
\end{align}
Letting $ \S_k(\phi_k)_T = \int_T \S_k(\phi_k) \dd x / |T|$ and $(\bold{\Pi}_k\bold{w})_T = \int_T \bold{\Pi}_k\bold{w} \dd x / |T|$ and $\bold{\varphi}_T = \int_{T} \bold{\varphi}\dd x / |T|$ for $T\in \cT_k$, by a scaled trace theorem, Poincar\'{e} inequality and estimates \eqref{int_operator-cr_approx}-\eqref{int_operator-cr_stab}, we further obtain that for any $F\in \mathcal{F}_k(\Omega)$,
\begin{align*}
 {\rm I}_F  &\leq \|[\S_k(\phi_k)] - [\S_k(\phi_k)]_F \|_{\bold{L}^2(F)} \|\bold{\varphi}-\bold{\varphi}_F\|_{\mathbb{L}^2(F)}    \\
        & \leq \frac{1}{2} \sum_{T\in\mathcal{N}_k(F)} \|\S_k(\phi_k) - \S_k(\phi_k)_T\|_{\bold{L}^2(F)} \sum_{T\in\mathcal{N}_k(F)} \|\bold{\varphi} - \bold{\varphi}_T\|_{\mathbb{L}^2(F)} \nonumber \\
         & \leq c h_F \sum_{T\in \mathcal{N}_k(F)}\|\bold{\nabla}\S_k(\phi_k)\|_{\bold{L}^2(T)} \sum_{T\in \mathcal{N}_k(F)} \|\bold{\nabla}\bold{\varphi}\|_{\mathbb{L}^2(T)}.
\end{align*}
Likewise, we can derive for any $F\in \mathcal{F}_k\setminus \mathcal{F}_k(\Omega)$ and $F\subset \partial T$,
\begin{align*}
     {\rm III}_F    & \leq
         \left( \| \S_k(\phi_k) - \S_k(\phi_k)_T  \|_{\bold{L}^2(F)} + \| \bold{\Pi}_k \bold{w} - (\bold{\Pi}_k \bold{w})_T \|_{\bold{L}^2(F)} \right)
         \| \bold{\varphi} - \bold{\varphi}_T \|_{\mathbb{L}^2(F)} \\
         & \leq c h_F ( \| \bold{\nabla} \S_k(\phi_k) \|_{\bold{L}^2(T)} +  \|\bold{\nabla} \bold{\Pi}_k\bold{w}\|_{\bold{L}^2(T)} ) \|\bold{\nabla} \bold{\varphi}\|_{\mathbb{L}^2(T)},\\
       {\rm II}_F   & \leq c h_F   \|\bold{\nabla} \bold{w}\|_{\bold{L}^2(T)} \|\bold{\nabla} \bold{\varphi}\|_{\mathbb{L}^2(T)}.
        \end{align*}
Since $\{\|\S_k(\phi_k)\|_{1,k}\}_{k\geq 0}$ and $\{\|\bold{\Pi}_k\bold{w}\|_{1,k}\}_{k\geq 0}$ are uniformly bounded, the last three estimates, the local quasi-uniformity of $\cT_k$ and $ \lim_{k\to\infty}\|h_k\|_{L^\infty(\Omega)}= 0$ imply
 \begin{align*}
&\lim_{k\to\infty}  \bigg| \sum_{F\in\mathcal{F}_{k}(\Omega)} {\rm I}_F \bigg| \leq \lim_{k\to\infty} c \|h_k\|_{L^\infty(\Omega)} \|\bold{\nabla}\bold{\varphi}\|_{\mathbb{L}^2(\Omega)} = 0\quad \mbox{and}\quad \lim_{k\to\infty} \sum_{F\in\mathcal{F}_k\setminus\mathcal{F}_{k}(\Omega)}{\rm III}_F + \sum_{F\in\mathcal{F}_k\setminus\mathcal{F}_{k}(\Omega)}{\rm II}_F  = 0,
\end{align*}
    which, together with \eqref{lem:vel_weak-conv_pf03}, yields
    \begin{equation}\label{lem:vel_weak-conv_pf04}
     \lim_{k\to\infty}    \sum_{F\in\mathcal{F}_{k}} ([\S_k(\phi_k)],\bold{\varphi}\cdot\bold{n}_F)_{L^2(F)}  = (\bold{w}, \bold{\varphi} \cdot \bold{n})_{L^2(\partial\Omega)}, \quad \forall \bold{\varphi} \in C^\infty(\overline{\Omega})^{d\times d}.
    \end{equation}
    Now passing to the limit $k\to\infty$ on both sides of \eqref{lem:vel_weak-conv_pf02} and noting \eqref{lem:vel_weak-conv_pf01} and \eqref{lem:vel_weak-conv_pf04} lead to
    \begin{equation}\label{lem:vel_weak-conv_pf05}
        (\overline{\bold{u}}, \mathbf{div}\,\bold{\varphi}) = - (\bold{\sigma},\bold{\varphi}) + (\bold{w}, \bold{\varphi}\cdot \bold{n} )_{L^2(\partial\Omega)}, \quad \forall \bold{\varphi} \in C^\infty(\overline{\Omega})^{d\times d}.
    \end{equation}
    Then, with $\bold{\varphi} \in C_0^\infty(\Omega)^{d\times d}$, $\bold{\sigma} = \bold{\nabla} \overline{\bold{u}}$ and $\overline{\bold{u}} \in \bold{H}^1(\Omega)$.

    \noindent\textit{Step 2}. $\overline{\bold{u}} = \bold{g}$ on $\partial \Omega$. Let $\bold{H}(\mathrm{div}):=\{\boldsymbol{v}\in
    \bold{L}^{2}(\Omega)~|~\mathrm{div}\,\boldsymbol{v}\in L^{2}(\Omega)\}$ and $\mathbb{H}(\mathbf{div}) :=\{\bold{\tau}\in \mathbb{L}^2(\Omega)~|~ \mathbf{div}\,\bold{\tau}\in \bold{L}^2(\Omega)\}$. We invoke the density of $C^\infty(\overline{\Omega})^d$ in $\bold{H}(\mathrm{div})$ and the unique {continuous} extension of the normal trace mapping $\gamma_n: C^\infty(\overline{\Omega})^d \to H^{-1/2}(\partial\Omega) $ with $\bold{v} \mapsto \bold{v} \cdot \bold{n} $ to $\gamma_n: \bold{H}(\mathrm{div}) \to H^{-1/2}(\partial\Omega) $ \cite[Theorems I.2.4 and I.2.5, p. 27]{GiraultRaviart:1986}. Applying these two results row-wisely to each $\bold{\varphi} \in C^\infty(\overline{\Omega})^{d\times d}$ in \eqref{lem:vel_weak-conv_pf05} gives
  $$      (\overline{\bold{u}}, \mathbf{div}\,\bold{\varphi}) = - (\bold{\nabla} \overline{\bold{u}}, \bold{\varphi} ) +  \langle \bold{\varphi}\cdot\bold{n}, \bold{w}\rangle_{L^2(\partial\Omega)} , \quad \forall \bold{\varphi} \in \mathbb{H}(\mathbf{div}),
  $$
    where $\langle\cdot,\cdot\rangle_{L^2(\partial\Omega)}$ is understood as the duality between $H^{-1/2}(\partial\Omega)^d$ and $H^{1/2}(\partial\Omega)^d$.
    Meanwhile, applying \cite[Theorems I.2.4 and I.2.5]{GiraultRaviart:1986} row-wisely also leads to Green's formula
        $(\overline{\bold{u}}, \mathbf{div}\, \bold{\varphi}) + (\bold{\nabla} \overline{\bold{u}}, \bold{\varphi}) = \langle \bold{\varphi}\cdot\bold{n}, \overline{\bold{u}}\rangle_{L^2(\partial\Omega)}$ for $\bold{\varphi}\in \mathbb{H}(\mathbf{div})$.
Consequently, $\langle \bold{\varphi}\cdot\bold{n}, \overline{\bold{u}}\rangle_{L^2(\partial\Omega)}
     = \langle \bold{\varphi}\cdot\bold{n}, \bold{w}\rangle_{L^2(\partial\Omega)} = \langle \bold{\varphi}\cdot\bold{n}, \bold{g}\rangle_{L^2(\partial\Omega)}$ for any $\bold{\varphi}\in \mathbb{H}(\mathbf{div}).$ This and  the surjectivity of $\gamma_n: \bold{H}(\mathrm{div}) \to H^{-1/2}(\partial\Omega)$ \cite[Corollary I.2.8, p. 28]{GiraultRaviart:1986} give $\overline{\bold{u}}=\bold{g}$ on $\partial\Omega$.

     \noindent\textit{Step 3}. $\mathrm{div}\, \overline{\bold{u}} = 0$. {By arguing as in  \eqref{lem:vel_weak-conv_pf02}, 
     we have
     \begin{align}\label{lem:vel_weak-conv_pf08}
        (\S_k(\phi_k), \bold{\nabla} \varphi) = - \sum_{T\in\cT_k}  (\mathrm{div}\, \S_k(\phi_k), \varphi)_{L^2(T)} + \sum_{F\in \mathcal{F}_k(\Omega)} ([\S_k(\phi_k)]\cdot \bold{n}_F, \varphi)_{L^2(F)}  , \quad \forall \varphi \in C^\infty_0(\Omega).
\end{align}}
Inserting $\varphi_F = \int_F \varphi \dd s / |F|$ and the mean value $[\S_k(\phi_k)]_F$ of $[\S_k(\phi_k)]$ over $F\in \mathcal{F}_k(\Omega)$, we derive similarly
\begin{align}
        & \quad \sum_{F\in \mathcal{F}_k(\Omega)} ( [\S_k(\phi_k)]\cdot \bold{n}_F, \varphi)_{L^2(F)}  = \sum_{F\in \mathcal{F}_k(\Omega)} ([\S_k(\phi_k)]\cdot \bold{n}_F, \varphi -  \varphi_F)_{L^2(F)}   \nonumber \\
        & = \sum_{F\in \mathcal{F}_k(\Omega)} ( [\S_k(\phi_k)] - [\S_k(\phi_k)]_F)\cdot \bold{n}_F , \varphi -  \varphi_F )_{L^2(F)}\to 0\quad \text{as}~k\to \infty. \label{lem:vel_weak-conv_pf09}
        \end{align}
Since $\mathrm{div}_k\, \S_k(\phi_k) = 0$, passing to the limit $k\to\infty$ in \eqref{lem:vel_weak-conv_pf08} and {noting the second weak convergence in \eqref{lem:vel_weak-conv_pf01} and \eqref{lem:vel_weak-conv_pf09}}, we obtain {$(\overline{\bold{u}}, \bold{\nabla} \varphi) = 0$} for any $\varphi \in C_0^\infty(\Omega)$. By the definition of the distributional divergence, we have $\mathrm{div}\,\overline{\bold{u}} = 0$.

\noindent \textit{Step 4.}  $\overline{\bold{u}} = \S(\phi)$. We employ the density of $C^\infty_0(\Omega)^d \cap \bold{V}$ in $\bold{V}$ \cite[Corollary I.2.5, p. 26]{GiraultRaviart:1986}. For any $\bold{v} \in C^\infty_0(\Omega)^d \cap \bold{V}$, $\bold{\Pi}_k\bold{v}\in \bold{V}_k$ due to $\mathrm{div}\bold{v} = 0$ and \eqref{int_operator-cr} over $\cT_k$. Standard interpolation error estimates  \cite{Brenner:2003,Ciarlet:2002,CrouzeixRaviart:1973} imply
    \begin{equation}\label{lem:vel_weak-conv_pf10}
        \| \bold{v} - \bold{\Pi}_k \bold{v} \|_{\bold{L}^2(\Omega)} + \|h_k\|_{L^\infty(\Omega)} \|\bold{v} - \bold{\Pi}_k \bold{v}\|_{1,k} \leq c \|h_k\|^2_{L^\infty(\Omega)} |\bold{v}|_{\bold{H}^2(\Omega)}.
    \end{equation}
     We may insert $\bold{v}_k=\bold{\Pi}_k\bold{v}$ in \eqref{disvp_stokes_div-free} over $\cT_k$, i.e.,
    \begin{equation}\label{lem:vel_weak-conv_pf11}
        \mu (\bold{\nabla}_k \S_k(\phi_k), \bold{\nabla}_k \bold{\Pi}_k \bold{v}) + (\alpha_\eps(\phi_k) \S_k(\phi_k), \bold{\Pi}_k \bold{v}) = (\bold{f}, \bold{\Pi}_k\bold{v}).
    \end{equation}
    With $k\to \infty$ in \eqref{lem:vel_weak-conv_pf11}, the weak convergence $\bold{\nabla}_k\S_k(\phi_k) \rightarrow \bold{\nabla} \overline{\bold{u}}$ in $\mathbb{L}^2(\Omega)$ and strong convergence $\lim_{k\to\infty}\| \bold{v} - \bold{\Pi}_k \bold{v} \|_{\bold{L}^2(\Omega)} + \|\bold{v} - \bold{\Pi}_k \bold{v}\|_{1,k} = 0$ (due to \eqref{lem:vel_weak-conv_pf10} and $\|h_k\|_{L^\infty(\Omega)} \to 0$) ensure
    \begin{equation}\label{lem:vel_weak-conv_pf12}
        \lim_{k\to\infty}\mu (\bold{\nabla}_k \S_k(\phi_k), \bold{\nabla}_k \bold{\Pi}_k \bold{v}) = \mu (\bold{\nabla} \overline{\bold{u}}, \bold{\nabla} \bold{v})\quad \mbox{and} \quad
        \lim_{k\to\infty}(\bold{f},\bold{\Pi}_k \bold{v}) = (\bold{f}, \bold{v}).
    \end{equation}
    In view of the pointwise convergence of $\{\phi_k\}_{k\geq 0}$ in \eqref{lem:vel_weak-conv_pf01}, Assumption \ref{ass:problem}(i) on $\alpha_\eps$ and  Lebesgue dominated convergence theorem, $\alpha_\eps (\phi_k) \bold{v} \to \alpha_\eps (\phi) \bold{v}$ in $\bold{L}^2(\Omega)$, which, together with $\S_k(\phi_k)\to \overline{\bold{u}}$ weakly in $\bold{L}^2(\Omega)$, further yields
    \[
       \lim_{k\to\infty} (\alpha_\eps (\phi_k) \S_k(\phi_k), \bold{v}) = ( \alpha_\eps(\phi) \overline{\bold{u}}, \bold{v}),\quad \forall \bold{v}\in C^\infty_0(\Omega)^d \cap \bold{V}.
    \]
    The above argument also shows that $\alpha_\eps (\phi_k) \S_k(\phi_k) \rightarrow \alpha_\eps (\phi) \overline{\bold{u}}$ weakly in $\bold{L}^2(\Omega).$
    This and the convergence $\lim_{k\to\infty}\|\bold{v} - \bold{\Pi}_k \bold{v}\|_{\bold{L}^2(\Omega)}=0$ yield
    \begin{equation}\label{lem:vel_weak-conv_pf13}
    \lim_{k\to\infty}    (\alpha_\eps(\phi_k) \S_k(\phi_k), \bold{\Pi}_k \bold{v}) = (\alpha_\eps(\phi) \overline{\bold{u}}, \bold{v}).
    \end{equation}
    Consequently, it follows from \eqref{lem:vel_weak-conv_pf12} and \eqref{lem:vel_weak-conv_pf13} that
    \[
        \mu(\bold{\nabla} \overline{\bold{u}}, \bold{\nabla} \bold{v}) + (\alpha_\eps(\phi) \overline{\bold{u}}, \bold{v})  = (\bold{f}, \bold{v}), \quad \forall \bold{v} \in C_0^\infty(\Omega)^d\cap\bold{V}.
    \]
     Since $C_0^\infty(\Omega)^d\cap\bold{V}$ is dense in $\bold{V}$, we arrive at $\overline{\bold{u}} = \S(\phi)$. Finally, the uniqueness of $\S(\phi)$ and a standard subsequence contradiction argument imply that  \eqref{vel_weak-conv} holds for the whole sequence.
\end{proof}

The {subsequent} argument follows the proof for the existence result of problem \eqref{min_phase-field}-\eqref{vp_stokes_div-free} \cite{GarckeHecht:2016b}, in which Sobolev compact embedding theorem plays an essential role. However, the compactness does not hold automatically for the CR FE space. To overcome the challenge, we define a vectorial version of the connection operator between the scalar CR FE space and the scalar conforming $P_d$ (quadratic for $d=2$ or cubic for $d=3$) FE space \cite{Brenner:2003}. Let
$$\bold{X}_k^c: = \{\bold{v}\in
\bold{H}^{1}(\Omega)~|~\bold{v}|_{T}\in P_{d}(T)^d~\forall T\in\mathcal{T}_k\}$$
and $\E_k : \bold{X}_k\mapsto \bold{X}_k^c$ be defined by
\begin{equation}\label{def:con_op}
    \E_{k} \bold{v}(p)=\dfrac{1}{\#\mathcal{N}_{p}}\displaystyle{\sum_{T\in\mathcal{N}_{p}}}\bold{v}|_{T}(p),\quad\hbox{for node $p$},
\end{equation}
where $\mathcal{N}_{p}$ is the set of $T\in\mathcal{T}_k$ that share a common node $p$ and $\#\mathcal{N}_p$ is the cardinality of $\mathcal{N}_{p}$. Note that $\E_k\bold{v}(p)=\bold{v}(p)$ at all centers $p$ of faces/edges in $\mathcal{F}_k$ since $\bold{v}\in \bold{X}_k$ is continuous at these nodes. $\E_k$ satisfies the following estimates
\cite[(3.3) in Lemma 3.2 and (3.10)-(3.11) in Corollary 3.3]{Brenner:2003}
\begin{equation}\label{con_operator-stab}
    \| \bold{v} - \E_k \bold{v}  \|_{\bold{L}^2(\Omega)} \leq c \|h_k\|_{L^\infty(\Omega)} \|\bold{v}\|_{1,k}, \quad
    \|\E_k \bold{v}\|_{\bold{L}^2(\Omega)} \leq c \|\bold{v}\|_{\bold{L}^2(\Omega)},\quad |\E_k \bold{v}|_{\bold{H}^1(\Omega)} \leq c \|\bold{\nabla}_k\bold{v}\|_{\bold{L}^2(\Omega)},\quad    \forall \bold{v} \in \bold{X}_k,
\end{equation}
where the constant $c$ depends only on the shape regularity of $\cT_k$.

\begin{lemma}\label{lem:vel_L2-strong-conv}
Let Assumption \ref{ass:problem}(i) and (iv) hold. Let $\{\phi_k\in \U_k\}_{k\geq 0}$ converge to some $\phi \in \U $ weakly in $H^1(\Omega)$ as $\lim_{k\to\infty}\|h_k\|_{L^\infty(\Omega)}= 0$. Then the sequence of discrete velocities $\{\S_k(\phi_k) \in \bold{\Pi}_k \bold{w} + \bold{V}_k\}_{k\geq0}$ converges strongly in $\bold{L}^2(\Omega)$ to the unique solution  $\S(\phi) \in \bold{U}$ of problem \eqref{vp_stokes_div-free}.
\end{lemma}
\begin{proof}
Applying the operator $\E_k$ to $\bold{u}_k = \S_k(\phi_k)$ and using the stability estimates in \eqref{con_operator-stab} lead to
\[      \|\E_k \bold{u}_k\|_{\bold{H}^1(\Omega)} \leq c (\|\bold{u}_k\|_{\bold{L}^2(\Omega)} + \|\bold{\nabla}_k\bold{u}_k\|_{\bold{L}^2(\Omega)}).
\]
Since $\{\|\bold{u}_k\|_{\bold{L}^2(\Omega)} + \|\bold{u}_k\|_{1,k}\}_{k\geq 0}$ is uniformly bounded (cf. the proof of Lemma \ref{lem:vel_weak-conv}), by Sobolev compact embedding theorem, we can extract a subsequence, still denoted by $\{\E_k \bold{u}_k\}_{k\geq 0}$, and find $\bold{u} \in \bold{H}^1(\Omega)$ such that
$  \lim_{k\to\infty}\| \E_k \bold{u}_k -\bold{u}\|_{\bold{L}^2(\Omega)}=0.$
It follows from the error estimate in \eqref{con_operator-stab} and the uniform boundedness of $\{\|\bold{u}_k\|_{1,k}\}_{k\geq 0}$ again that
$    \lim_{k\to\infty}    \|\bold{u}_k - \E_k \bold{u}_k\|_{\bold{L}^2(\Omega)}=0.$
These two estimates and the triangle inequality imply $\S_k(\phi_k) = \bold{u}_k \to \bold{u} $ in $\bold{L}^2(\Omega)$. By noting the $\bold{L}^2(\Omega)$ weak convergence of $\{\S_k(\phi_k)\}_{k\geq0}$ to $\S(\phi)$ in Lemma \ref{lem:vel_weak-conv} and the uniqueness of the limit, we get the desired assertion for the whole sequence with $\bold{u} = \S(\phi)$.
\end{proof}

\begin{remark}
The proof retrieves the discrete compactness of CR elements \cite{Stummel:1980} that every bounded sequence $\{\|\bold{v}_k\|_{\bold{L}^2(\Omega)}+\|\bold{\nabla}_k\bold{v}_k\|_{\bold{L}^2(\Omega)}\}_{k\geq0}$ with each $\bold{v}_k \in \bold{X}_k$ has a subsequence converging in $\bold{L}^2(\Omega)$ to some $\bold{v} \in \bold{H}^1(\Omega)$ as $\|h_k\|_{L^\infty(\Omega)}\to 0$.
\end{remark}

\begin{lemma}\label{lem:vel_strong-conv}
Let Assumption \ref{ass:problem}(i) and (iv) hold. Let $\{\phi_k\in \U_k\}_{k\geq 0}$ converge to some $\phi \in \U $ weakly in $H^1(\Omega)$ as $\lim_{k\to\infty}\|h_k\|_{L^\infty(\Omega)}=0$. If $\|\bold{\nabla}_k(\bold{w} - \bold{\Pi}_k \bold{w})\|_{\bold{L}^2(\Omega)} \to 0$, the sequence of discrete velocities $\{\S_k(\phi_k) \in \bold{\Pi}_k \bold{w} + \bold{V}_k\}_{k\geq0}$ converges to the unique solution  $\S(\phi) \in \bold{U}$ of \eqref{vp_stokes_div-free} in the sense that $\lim_{k\to\infty}\|\bold{\nabla}_k(\S_k(\phi_k) - \S(\phi))\|_{\bold{L}^2(\Omega)}= 0$.
\end{lemma}

\begin{proof}
    In the identity \eqref{disvp_stokes_div-free} with $\phi_\cT = \phi_k$ over $\cT_k$, we take $\bold{v}_k = \S_k(\phi_k) - \bold{\Pi}_k\bold{w} $ and obtain
\begin{align}
        \mu \|\bold{\nabla}_k\S_k(\phi_k)\|_{\bold{L}^2(\Omega)}^2   + (\alpha_\eps(\phi_k),|\S_k(\phi_k)|^2) =& ( \bold{f}, \S_k(\phi_k) - \bold{\Pi}_k \bold{w}) \nonumber\\
         & + \mu (\bold{\nabla}_k \S_k(\phi_k), \bold{\nabla}_k \bold{\Pi}_k \bold{w} ) + ( \alpha_\eps(\phi_k) \S_k(\phi_k),\bold{\Pi}_k \bold{w}).\label{lem:vel_strong-conv_pf01}
        \end{align}
Due to \eqref{int_operator-cr_approx}, $\lim_{k\to\infty}\|\bold{w} - \bold{\Pi}_k \bold{w}\|_{\bold{L}^2(\Omega)}= 0$. From this, the assumption $\lim_{k\to\infty}\|\bold{\nabla}_k(\bold{w} - \bold{\Pi}_k \bold{w})\|_{\bold{L}^2(\Omega)} = 0$ and Lemma \ref{lem:vel_weak-conv}, it follows that
    \[
 \lim_{k\to\infty}       (\bold{f},\S_k(\phi_k) - \bold{\Pi}_k \bold{w}) + \mu (\bold{\nabla}_k\S_k(\phi_k), \bold{\nabla}_k \bold{\Pi}_k \bold{w} ) = (\bold{f}, \S(\phi) - \bold{w}) + \mu (\bold{\nabla}\S(\phi), \bold{\nabla} \bold{w}).
    \]
    By the argument for  \eqref{lem:vel_weak-conv_pf13} (which holds independent of the subsequence), we have
$ \lim_{k\to\infty}       (\alpha_\eps(\phi_k) \S_k(\phi_k), \bold{\Pi}_k \bold{w}) =  (\alpha_\eps(\phi) \S(\phi),   \bold{w})$.
Passing to the limits on both sides of \eqref{lem:vel_strong-conv_pf01} gives
\begin{equation}\label{lem:vel_strong-conv_pf02}
        \begin{aligned}
    \lim_{k\to\infty} \mu  \|\bold{\nabla}_k\S_k(\phi_k)\|_{\bold{L}^2(\Omega)}^2   + (\alpha_\eps(\phi_k), |\S_k(\phi_k)|^2) =
 \mu (\bold{f},\S(\phi) - \bold{w})  + \mu ( \bold{\nabla}\S(\phi),\bold{\nabla} \bold{w})
         + (\alpha_\eps(\phi) \S(\phi),\bold{w}).
        \end{aligned}
    \end{equation}
Next by noting that $\bold{u} = \S(\phi)$ solves problem \eqref{vp_stokes_div-free} and inserting $\bold{v} = \S(\phi) - \bold{w} \in \bold{V}$, we get
    \begin{equation}\label{lem:vel_strong-conv_pf03}
     \mu \|\bold{\nabla} \S(\phi)\|_{\bold{L}^2(\Omega)}^2 + ( \alpha_\eps(\phi),|\S(\phi)|^2)
      = \mu (\bold{f},\S(\phi) - \bold{w}) + \mu ( \bold{\nabla}\S(\phi),\bold{\nabla} \bold{w})
         + (\alpha_\eps(\phi) \S(\phi), \bold{w}).
    \end{equation}
    The last two identities imply
    \begin{equation}\label{lem:vel_strong-conv_pf04}
      \lim_{k\to\infty}  \mu \|\bold{\nabla}_k\S_k(\phi_k)\|_{\bold{L}^2(\Omega)}^2  + (\alpha_\eps(\phi_k), |\S_k(\phi_k)|^2) = \mu \|\bold{\nabla} \S(\phi)\|_{\bold{L}^2(\Omega)}^2  + ( \alpha_\eps(\phi),|\S(\phi)|^2).
    \end{equation}
    Further, using the pointwise convergence of $\{\phi_k\}_{k\geq0}$ to $\phi$ (up to a subsequence), Assumption \ref{ass:problem}(i) on $\alpha_\eps$,  $\bold{L}^2$ convergence of $\{\S_k(\phi_k)\}_{k\geq 0}$, Lebesgue dominated convergence theorem and the subsequence contradiction argument, we deduce that for the whole sequence,
        \begin{align*}
         &\lim_{k\to\infty}(\alpha_\eps(\phi_k),|\S_k(\phi_k)|^2) - ( \alpha_\eps(\phi), |\S(\phi)|^2)\\
        =& \lim_{k\to\infty}(\alpha_\eps(\phi_k), |\S_k(\phi_k)|^2 - |\S(\phi)|^2) + ((\alpha_\eps(\phi_k) - \alpha_\eps(\phi)), |\S(\phi)|^2) = 0,
        \end{align*}
    which, along with \eqref{lem:vel_strong-conv_pf04}, yields $\lim_{k\to\infty}\mu \|\bold{\nabla}_k\S_k(\phi_k)\|_{\bold{L}^2(\Omega)}^2 =  \mu \|\bold{\nabla}\S(\phi)\|_{\bold{L}^2(\Omega)}^2$. This and $\bold{\nabla}_k\S_k(\phi_k)\rightharpoonup \bold{\nabla}\S(\phi)$ weakly in $\mathbb{L}^2(\Omega)$ in Lemma \ref{lem:vel_weak-conv} complete the proof of the lemma.
\end{proof}

Now we can state the main result of this section.

\begin{theorem}\label{thm:conv_uniform}
Let Assumption \ref{ass:problem} hold. Let $\{\cT_k\}_{k\geq 0}$ be a sequence of uniformly refined meshes and $\{(\phi_k^\ast, \bold{u}_k^\ast)\}_{k\geq 0}$ a sequence of discrete minimizers to problem \eqref{dismin_phase-field}-\eqref{disvp_stokes_div-free} and the associated discrete velocity fields. If $\lim_{k\to\infty}\|\bold{\nabla}_k(\bold{w} - \bold{\Pi}_k\bold{w})\|_{\bold{L}^2(\Omega)}=0$, then there exists a subsequence of pairs $\{(\phi_{k_j}^\ast,\bold{u}_{k_j}^\ast)\}_{j\geq 0}$ converging strongly to a minimizer $\phi^\ast \in \U$ to problem \eqref{min_phase-field}-\eqref{vp_stokes_div-free} in $H^1(\Omega)$ and the associated velocity field $\bold{u}^\ast\in \bold{U}$ in the sense that
    \begin{equation}\label{vel_strong-conv}
       \lim_{j\to\infty} \big(\|\bold{u}^\ast_{k_j} - \bold{u}^\ast \|_{\bold{L}^2(\Omega)} + \|\bold{\nabla}_{k_j}(\bold{u}^\ast_{k_j} - \bold{u}^\ast) \|_{\bold{L}^2(\Omega)}\big) = 0.
    \end{equation}
\end{theorem}

\begin{proof}
    Since $\beta \in \U_k$ for each $k\geq 0$, we may insert $\phi_k = \beta$ in \eqref{disvp_stokes_div-free} over $\cT_k$ and deduce from the convergence in Lemmas \ref{lem:vel_L2-strong-conv} and \ref{lem:vel_strong-conv} and the continuity in Assumption \ref{ass:problem}(ii) that    \begin{equation}\label{thm:conv_uniform_pf01}
        \sup_{k\geq 0} \|\S_k(\beta)\|_{\bold{L}^2(\Omega)} + \sup_{k\geq0} \left | (G(\bold{x},\S_k(\beta), \bold{\nabla}_k\S_k(\beta)),1) \right | < \infty.
    \end{equation}
Meanwhile, $\int_\Omega G(\bold{x}, \bold{u}_k^\ast, \bold{\nabla}_k\bold{u}_k^\ast) \dx$ is bounded from below by Assumption \ref{ass:problem}(iii). Hence from \eqref{thm:conv_uniform_pf01}, Assumption \ref{ass:problem}(i) on $\alpha_\eps$, we deduce that there exist two constants $c_1,c_2>0$ such that
    \[
        -c_1 + \gamma \frac{\eps}{2} |\phi_k^\ast|^2_{H^1(\Omega)} \leq -c_1 + \gamma\mathcal{P}_\eps(\phi^\ast_k)  \leq \J_{k}^\eps(\phi_k^\ast) \leq \J_{k}^\eps(\beta) \leq c_2 + \tfrac{\gamma}{\eps}( f(\beta),1) , \quad \forall k\geq 0,
    \]
    which, along with the bound $\|\phi_k^\ast\|_{L^2(\Omega)} \leq \sqrt{|\Omega|}$, implies that the sequence $\{\phi_k^\ast\}_{k\geq0}\subset \U$ is uniformly bounded in $H^1(\Omega)$. Since the set $\U$ is closed and convex, the reflexivity of the space $H^1(\Omega)$ and Sobolev compact embedding theorem yield a subsequence $\{\phi_{k_j}^\ast\}_{j\geq 0}$ and some $\overline{\phi} \in \U$ such that
\begin{equation}\label{thm:conv_uniform_pf03}
        \phi_{k_j}^\ast \rightarrow \overline{\phi}\quad \text{weakly in}~H^1(\Omega),\quad\phi_{k_j}^\ast \to \overline{\phi}\quad \text{a.e. in}~\Omega.
    \end{equation}
    By Lemmas \ref{lem:vel_L2-strong-conv} and \ref{lem:vel_strong-conv} with $\bold{u}_{k_j}^\ast=\S_{k_j}(\phi^\ast_{k_j})$ and $\overline{\bold{u}}=\S(\overline{\phi})$ and  Assumption \ref{ass:problem}(ii), we obtain
    \begin{equation}\label{thm:conv_uniform_pf04}
     \lim_{j\to\infty}   (G(\bold{x}, \bold{u}_{k_j}^\ast, \bold{\nabla}_{k_j} \bold{u}_{k_j}^\ast),1) = ( G(\bold{x}, \overline{\bold{u}}, \bold{\nabla}\overline{\bold{u}}),1),
    \end{equation}
    while Lemma \ref{lem:vel_L2-strong-conv}, the pointwise convergence in \eqref{thm:conv_uniform_pf03} and Lebesgue dominated convergence theorem imply
\begin{equation}\label{thm:conv_uniform_pf05}
       \lim_{j\to\infty} (\alpha_{\eps}(\phi^\ast_{k_j}), |\bold{u}_{k_j}^\ast|^2)+ \tfrac{\gamma}{\eps}(f (\phi_{k_j}^\ast),1)\to
        (\alpha_{\eps}(\overline{\phi}), |\overline{\bold{u}}|^2) + \tfrac{\gamma}{\eps} (f (\overline{\phi}),1).
    \end{equation}
This argument yields a \textit{liminf} inequality, one of two main ingredients in $\Gamma$-convergence \cite{Braides:2002}. For a \textit{limsup} inequality, i.e., the construction of a recovery sequence, we define an auxiliary admissible set $\widetilde{\U}:= \{\phi \in C^\infty(\overline{\Omega})~|~\phi \in [0,1)~\text{a.e. in}~\Omega, \int_{\Om}\phi\dx < \beta |\Om| \}$. The density of $C^\infty(\overline{\Omega})$ in $H^1(\Omega)$ via the mollifier  \cite{AdamsFournier:2003} implies that $C^\infty(\overline{\Omega})\cap\U$ is dense in the set $\U$ in the $H^1(\Omega)$-norm. This and the identity $\lim_{\lambda\to1^-}\|\lambda \phi - \phi\|_{H^1(\Omega)}= 0 $ for any $\phi \in H^1(\Omega)$ imply that the set $\widetilde{\U}$ is dense in $\U$ in the $H^1(\Omega)$-norm. Let $S_k$ be the $H^1(\Omega)$-conforming linear FE space. By invoking the classical nodal interpolation operator $I_k: H^2(\Omega) \rightarrow S_k$ \cite{BrennerScott:2008,Ciarlet:2002}, we have
    \begin{align}\label{thm:conv_uniform_pf06}
\lim_{k\to\infty}      \|\phi - I_k \phi\|_{H^1(\Omega)} &\to 0,\quad \forall \phi \in \widetilde{\U},\\
        \| \phi - I_k \phi \|_{L^1(\Omega)} &\leq \sqrt{|\Omega|} \| \phi - I_k \phi \|_{L^2(\Omega)} \leq c \sqrt{|\Omega|} \|h_k\|_{L^\infty(\Omega)}^{2} |\phi|_{H^2(\Omega)},\nonumber
    \end{align}
    which further implies
    \begin{equation}\label{thm:conv_uniform_pf07}
        ( I_k \phi,1)\leq ( \phi,1) + c \sqrt{|\Omega|} \|h_k\|_{L^\infty(\Omega)}^{2} |\phi|_{H^2(\Omega)},\quad  \forall \phi \in \widetilde{\U}.
    \end{equation}
    Since $\phi \in \widetilde{\U}$, clearly $I_k \phi \in [0,1)$ a.e. in $\Omega$ and in view of \eqref{thm:conv_uniform_pf07}, when $k$ is sufficiently large, i.e., $\|h_k\|_{L^\infty(\Omega)}$ is sufficiently small, $\int_\Omega I_k \phi \dd \bold{x} < \beta |\Omega|$. Hence $I_k \phi \in \U_k$. Using \eqref{thm:conv_uniform_pf06}, Lemmas \ref{lem:vel_L2-strong-conv} and \ref{lem:vel_strong-conv}, the continuity in Assumption \ref{ass:problem}(ii) and a standard subsequence contradiction argument, we deduce as in \eqref{thm:conv_uniform_pf04}-\eqref{thm:conv_uniform_pf05} and find
\begin{equation}\label{thm:conv_uniform_pf08}
     \lim_{k\to\infty}   \J_{k}^\eps(I_k\phi) = \J^\eps(\phi),\quad \forall \phi \in \widetilde{\U}.
    \end{equation}
    Now combining the weak convergence in \eqref{thm:conv_uniform_pf03}, \eqref{thm:conv_uniform_pf04}-\eqref{thm:conv_uniform_pf05}, the weak lower semi-continuity of $H^1(\Omega)$ semi-norm and \eqref{thm:conv_uniform_pf08} and the minimizing property $\J_k^\eps(\phi_k^\ast)\leq \J_k^\eps(I_k \phi) $ for any $\phi \in \widetilde{\U}$ and any sufficiently large $k$ gives
    \begin{equation}\label{thm:conv_uniform_pf09}
        \J^\eps (\overline{\phi}) \leq \liminf_{j\to\infty} \J_{k_j}^\eps(\phi_{k_j}^\ast) \leq \limsup_{j\to\infty} \J_{k_j}^\eps(\phi_{k_j}^\ast)
        \leq \limsup_{k\to\infty} \J_{k}^\eps(\phi_{k}^\ast) \leq \limsup_{k\to\infty} \J_{k}^\eps(I_k\phi) = \J^\eps(\phi), \quad \forall \phi \in \widetilde{\U}.
    \end{equation}
By the density of $\widetilde{\U}$ in $\U$ and Lemma \ref{lem:sol-map_cont}, we get $\phi^{\ast} = \overline{\phi}$ in the weak convergence of \eqref{thm:conv_uniform_pf03} and also the strong convergence in Lemmas \ref{lem:vel_L2-strong-conv} and \ref{lem:vel_strong-conv} with $\bold{u}^\ast = \S(\phi^\ast)$. By the density of $\widetilde{\U}$ in $\U$ and Lemma \ref{lem:sol-map_cont} (applied to a minimizer $\phi^\ast\in \U$), the argument in \eqref{thm:conv_uniform_pf09} implies that $\lim_{k\to\infty} \J_{k}^\eps(\phi_{k}^\ast) = \J^\eps(\phi^\ast)$. This, the weak convergence in \eqref{thm:conv_uniform_pf03} and \eqref{thm:conv_uniform_pf04}-\eqref{thm:conv_uniform_pf05} again yield the desired strong convergence of $\{\phi_{k_j}^\ast\}_{j\geq0}$.
\end{proof}

In practical numerical treatment of \eqref{vp_stokes_div-free}, the pressure field $p$ is approximated by piecewise constants. This motivates the analysis of discrete pressure fields $\{p_k^\ast\in Q_k\}_{k\geq0}$ associated with minimizing pairs $\{(\phi_k^\ast,\bold{u}_k^\ast)\}_{k\geq 0}$ for \eqref{dismin_phase-field}-\eqref{disvp_stokes_div-free}. The first key ingredient in the analysis is the following well-known inf-sup condition \cite{GiraultRaviart:1986}
\begin{equation}\label{inf-sup}
    \sup_{\bold{0}\neq\bold{v} \in \bold{H}_0^1(\Omega)} \frac{( \mathrm{div}\, \bold{v}, q)}{|\bold{v}|_{\bold{H}^1(\Omega)}} \geq c \| q \|_{L^2(\Omega)}, \quad \forall q \in L_0^2(\Omega)
\end{equation}
with $c > 0$ depending only on $\Omega$. Since the minimizer $\phi^\ast\in \U$ and the associated velocity field $\bold{u}^\ast \in \bold{U}$ solve \eqref{vp_stokes_div-free}, a bounded linear functional in the polar set $\bold{V}^{0}:=\{\bold{l}\in \bold{H}^{-1}(\Omega) ~|~\langle \bold{l},\bold{v} \rangle = 0~\forall \bold{v}\in \bold{V}\}$ is given by $(\phi^\ast, \bold{u}^\ast) \in \U \times \bold{U}$ through \eqref{vp_stokes_div-free}. Then it follows from \eqref{inf-sup} and Lemma I.4.1 in \cite{GiraultRaviart:1986} that there exists a unique $p^\ast \in L_0^2(\Omega)$ such that
\begin{subequations}
    \begin{align}
    \mu (\bold{\nabla} \bold{u}^\ast,\bold{\nabla} \bold{v}) + (\alpha_\eps(\phi^\ast) \bold{u}^\ast, \bold{v}) - ( \mathrm{div}\,\bold{v}, p^\ast)&= (\bold{f}, \bold{v}), \quad \forall \bold{v}  \in \bold{H}_0^1(\Omega), \label{vp_stokes_motion}\\
    \label{vp_stokes_incompress}
       (\mathrm{div}\,\bold{u}^\ast ,q ) &= 0, \quad \forall q \in L_0^2(\Omega).
    \end{align}
\end{subequations}
In parallel with \eqref{inf-sup}, there also holds a discrete version \cite{CrouzeixRaviart:1973}
\begin{equation}\label{inf-sup_disc}
    \sup_{\bold{0}\neq\bold{v}_k\in \bold{Z}_k} \frac{( \mathrm{div}_k \bold{v}_k, q_k)}{\|\bold{v}_k\|_{1,k}} \geq c \|q_k\|_{L^2(\Omega)},
\end{equation}
with $c$ depending only on the shape regularity of $\cT_k$. Condition \eqref{inf-sup_disc} ensures the existence of a unique $p_k^\ast \in Q_k$ associated with the minimizing pair $(\phi_k^\ast,\bold{u}_k^\ast)\in \U_k \times \bold{U}_k$ of problem \eqref{dismin_phase-field}-\eqref{disvp_stokes_div-free} such that
\begin{subequations}\label{disvp_stokes}
    \begin{align}\label{disvp_stokes_motion}
    \mu ( \bold{\nabla} \bold{u}_k^\ast , \bold{\nabla} \bold{v}_k) + (\alpha_\eps(\phi^\ast) \bold{u}_k^\ast , \bold{v}_k) - (\mathrm{div}_k \bold{v}_k, p_k^\ast) &= ( \bold{f}, \bold{v}_k), \quad \forall \bold{v}_k  \in \bold{Z}_k,\\
   \label{disvp_stokes_incompress}
        (\mathrm{div}_k\bold{u}^\ast_k, q_k) &= 0 ,\quad \forall q_k \in Q_k.
    \end{align}
\end{subequations}

We have the following convergence of discrete pressure fields $\{p_k^\ast\}_{k\geq0}\subset \bigcup_{k\geq 0}Q_k$ associated with minimizing pairs $\{(\phi_k^\ast,\bold{u}_k^\ast)\}_{k\geq 0}$ for problem \eqref{dismin_phase-field}-\eqref{disvp_stokes_div-free}.

\begin{theorem}\label{thm:conv_pre_uniform}
     Let $\{p_k^\ast\}_{k \geq 0}$ be the sequence of discrete pressure fields with each $p_k^\ast \in Q_k$ given by \eqref{disvp_stokes_motion}. Then under the assumption of Theorem \ref{thm:conv_uniform}, the subsequence $\{p_{k_j}^\ast\}_{j\geq 0}$ associated with the convergent subsequence $\{(\phi_{k_j}^\ast,\bold{u}_{k_j}^\ast)\}_{j\geq0}$ converges to $p^\ast\in L_0^2(\Omega)$ in \eqref{vp_stokes_motion} in $L^2(\Omega)$.
\end{theorem}
\begin{proof}
By taking $q = p^\ast - p_{k_j}^\ast$ in \eqref{inf-sup} and using \eqref{vp_stokes_motion}, we get
\begin{align}
    c\|p^\ast - p^\ast_{k_j}\|_{L^2(\Omega)} & \leq \sup_{\bold{0}\neq \bold{v} \in \bold{H}_0^1(\Omega)} \frac{(\mathrm{div}\,\bold{v},p^\ast - p^\ast_{k_j}) }{|\bold{v}|_{\bold{H}^1(\Omega)}} \nonumber \\
    & = \sup_{\bold{0}\neq \bold{v} \in \bold{H}_0^1(\Omega)} \frac{ \mu (\bold{\nabla} \bold{u}^\ast, \bold{\nabla} \bold{v}) + (\alpha_\eps(\phi^\ast) \bold{u}^\ast, \bold{v}) - ( \mathrm{div}\,\bold{v}, p_{k_j}^\ast) - (\bold{f}, \bold{v})}{|\bold{v}|_{\bold{H}^1(\Omega)}}.\label{thm:pre_strong_pf01}
\end{align}
Next we split the numerator in \eqref{thm:conv_pre_uniform} into
    \begin{align*}
     [\mu (\bold{\nabla}_k  ( \bold{u}^\ast - \bold{u}^\ast_{k_j} ), \bold{\nabla} \bold{v}) +  &((\alpha_\eps(\phi^\ast) \bold{u}^\ast - \alpha_\eps(\phi_{k_j}^\ast)\bold{u}_{k_j}^\ast), \bold{v})]\\
    & + [(\bold{\nabla}_k \bold{u}^\ast_{k_j}, \bold{\nabla} \bold{v}) + (\alpha_\eps(\phi_{k_j}^\ast) \bold{u}_{k_j}^\ast, \bold{v}) - (\mathrm{div}\, \bold{v} , p_{k_j}^\ast) - (\bold{f}, \bold{v})]:={\rm I} + {\rm II}.
    \end{align*}
Using $\bold{\Pi}_k$, \textcolor{red}{\eqref{disvp_stokes_motion}}, and elementwise integration by parts, and noting $(\mu\bold{\nabla}\bold{u}^\ast_k - p_{k_j}^\ast \bold{I}_d ) \cdot \bold{n}_{\partial T}$ ($\bold{I}_d$ is the $d \times d$ identity matrix) is a constant vector on $F\subset \partial T$, \eqref{int_operator-cr} and \eqref{int_operator-cr_approx}, we further get for any $\bold{v} \in \bold{H}_0^1(\Omega)$,
\begin{align}\label{thm:pre_strong_pf02}
  |{\rm II}|  & = \left|(\bold{\nabla}_k \bold{u}^\ast_{k_j}  , \bold{\nabla}_k ( \bold{v} - \bold{\Pi}_k \bold{v} )) + (\alpha_\eps(\phi_{k_j}^\ast) \bold{u}_{k_j}^\ast, \bold{v} - \bold{\Pi}_{k_j} \bold{v} ) - ( \mathrm{div}_k ( \bold{v} - \bold{\Pi}_{k_j} \bold{v} ), p_{k_j}^\ast) - (\bold{f},\bold{v} - \bold{\Pi}_{k_j} \bold{v} ) \right|\nonumber \\
    & = \left| \sum_{T\in\cT_{k_j}}\left( (\alpha_\eps(\phi_{k_j}^\ast) \bold{u}_{k_j}^\ast - \bold{f}) , \bold{v} - \bold{\Pi}_{k_j} \bold{v} )_{L^2(T)} + ( \bold{v} - \bold{\Pi}_{k_j} \bold{v} , (\mu\bold{\nabla}\bold{u}^\ast_{k_j} - p_{k_j}^\ast\bold{I}_d ) \cdot \bold{n}_{\partial T})_{L^2(\partial T)} \right) \right| \nonumber \\
    & \leq c \sum_{T\in \cT_k}h_T\|\alpha_\eps(\phi_k^\ast) \bold{u}_k^\ast - \bold{f}\|_{\bold{L}^2(T)}\|\bold{\nabla}\bold{v}\|_{\bold{L}^2(T)} \leq c \left(\sum_{T\in \cT_k}h_T^2\|\alpha_\eps(\phi_{k_j}^\ast) \bold{u}_{k_j}^\ast - \bold{f}\|_{\bold{L}^2(T)}^2\right)^{1/2} |\bold{v}|_{\bold{H}^1(\Omega)}.
\end{align}
Meanwhile, by Poincar\'{e} inequality, we deduce
\begin{equation}\label{thm:pre_strong_pf03}
     |{\rm I}|  \leq c \big(\| \bold{u}^\ast - \bold{u}_{k_j}^\ast \|_{1,k_j} + \| \alpha_\eps(\phi^\ast) \bold{u}^\ast - \alpha_\eps(\phi_{k_j}^\ast)\bold{u}_{k_j}^\ast \|_{\bold{L}^2(\Omega)}\big) |\bold{v}|_{\bold{H}^1(\Omega)},\quad \forall \bold{v} \in \bold{H}_0^1(\Omega).
\end{equation}
The preceding estimates give
\begin{equation}\label{thm:pre_strong_pf04}
    \|p^\ast - p^\ast_{k_j}\|_{L^2(\Omega)} \leq c \bigg(\Big(\sum_{T\in \cT_k}h_T^2\|\alpha_\eps(\phi_{k_j}^\ast) \bold{u}_{k_j}^\ast - \bold{f}\|_{\bold{L}^2(T)}^2\Big)^{1/2} + \| \bold{u}^\ast - \bold{u}_{k_j}^\ast \|_{1,k_j} + \| \alpha_\eps(\phi^\ast) \bold{u}^\ast - \alpha_\eps(\phi_{k_j}^\ast)\bold{u}_{k_j}^\ast \|_{\bold{L}^2(\Omega)} \bigg).
\end{equation}
Since Theorem \ref{thm:conv_uniform} implies that $\{ \|\bold{u}_{k_j}^\ast\|_{L^2(\Omega)}\}_{k\geq 0}$ is uniformly bounded, it follows from Assumption (i) that
\begin{equation}\label{thm:pre_strong_pf05}
\begin{aligned}
 \lim_{j\to\infty} \sum_{T\in \cT_{k_j}}h_T^2\|\alpha_\eps(\phi_{k_j}^\ast) \bold{u}_{k_j}^\ast - \bold{f}\|_{\bold{L}^2(T)}^2  &\leq \lim_{j\to\infty} \|h_{k_j}\|_{L^\infty(\Omega)}^2 \|\alpha_\eps(\phi_{k_j}^\ast) \bold{u}_{k_j}^\ast - \bold{f}\|_{\bold{L}^2(\Omega)}^2 \\
 &\leq \lim_{j\to\infty} c \|h_{k_j}\|_{L^\infty(\Omega)}^2 = 0.
\end{aligned}
\end{equation}
Since $\phi_{k_j}^\ast$ converges to $\phi^\ast$ weakly in $H^1(\Omega)$, Sobolev compact embedding theorem yields a pointwise convergent subsequence, still denoted by $\{\phi_{k_j}^\ast\}_{j\geq0}$. Then Lebesgue dominated convergence theorem and \eqref{vel_strong-conv} again imply
\[
 \lim_{j\to\infty}  \|\alpha_\eps(\phi^\ast) \bold{u}^\ast - \alpha_\eps(\phi_{k_j}^\ast)\bold{u}_{k_j}^\ast\|_{\bold{L}^2(\Omega)} \leq \lim_{j\to\infty}\| ( \alpha_\eps(\phi^\ast) - \alpha_\eps(\phi_{k_j}^\ast) )\bold{u}^\ast  \|_{\bold{L}^2(\Omega)} + \|\alpha_\eps(\phi_{k_j}^\ast) (\bold{u}^\ast - \bold{u}_{k_j}^\ast)\|_{\bold{L}^2(\Omega)} = 0,
\]
which also holds for the original subsequence by a standard contradiction argument. Collecting this, \eqref{vel_strong-conv}, \eqref{thm:pre_strong_pf04} and \eqref{thm:pre_strong_pf05} gives  the desired assertion.
\end{proof}

\section{Numerical results and discussions}

Now we illustrate the nonconforming CR-$P_{0}$ elements and the conforming $P_2$-$P_1$ elements on several examples to demonstrate the advantage of the proposed approach. The objective functional $G$ is given by the dissipated energy $\frac{1}{2} \int_\Omega|\bold{\nabla} \bold{u}|^2 \dx$  with $\mu=1$ and $\bold{f}=\bold{0}$ in \eqref{vp_stokes_div-free}.  $f(\phi)$ is a double-well potential, cf. section \ref{subsect:solver}, and $\alpha_\eps(\phi) = 10000 (1-\phi)^2$. The theoretical analysis focuses on the pure Dirichlet boundary condition $\bold{u}=\bold{g}$. Since practical applications often involve a mixed-type boundary condition, we prescribe the velocity of the fluid as the inlet boundary data and impose a Neumann boundary condition ($(\mu\bold{\nabla u}-p \bold{I}_d)\cdot\bold{n} = \bold{0}$) at the outlet while the remaining part is equipped with a zero Dirichlet boundary condition in the numerical examples. All numerical simulations are performed using MATLAB R2023a on a standard desktop with a 13th Gen Intel(R) Core(TM) i7-13700 2.10 GHz CPU and 32GB memory. The PDE Toolbox and the MATLAB package iFEM \cite{Chen:2009} are used in the experiments. Consider the following five cases in two dimensions (2d) and two cases in three dimensions (3d):
\begin{enumerate}[label=(\alph*)]
  \item (\textit{Pipe Bend}) The domain $\Omega$ is a square $(0,1)^2$. The inlet boundary condition is $\bold{g}=(1,0)^{\mathrm{T}}$ on $\{(x,y)~|~x=0, 0.7\leq y \leq 0.9\}$. A zero Neumann boundary condition is specified on $\{(x,y)~ |~ 0.7\leq x \leq 0.9, y=0\}$, and a zero Dirichlet boundary condition on the remaining part of the boundary $\partial \Omega$. \label{exp:pipebend}

  \item (\textit{Left Inflow}) The setting is identical with case \ref{exp:pipebend}, but the inlet boundary condition is $\bold{g}=\left(4y\left(1-y\right),0\right)^{\mathrm{T}}$ on the left side $x=0$, and the outlet boundary lies at $\{(x,y)| x=1, 0.3\leq y \leq 0.7\}$. \label{exp:leftinflow}

  \item (\textit{Three Inflows}) The domain $\Omega$ is a square $(0,1)^2$. The three inflows are positioned in the middle of the top, bottom, and left edges, each spanning a width of 0.2. The flow directions are downward, upward and rightward, respectively, with a unit speed. The outflow is located in the middle of the right edge with a spanning width of 0.2. The remaining boundaries are subject to a zero Dirichlet boundary condition. \label{exp:threeinflows}

  \item (\textit{Rugby}) The domain is a rectangle $(-0.5,1.5)\times(-0.5,0.5)$. A Dirichlet boundary $\bold{g}=\left(-\left(y-0.5\right)\left(y+0.5\right),0\right)^{\mathrm{T}}$ is applied at $x=-0.5$, and a zero Neumann boundary condition on the whole right side. \label{exp:rugby}

  \item (\textit{Bypass}) The domain $\Omega$ is a rectangle $(0,1.5) \times(-0.5,0.5)$. Two inlets are located at $\{(x,y)~|~x=0, 0.15 \leq y  \leq 0.35\} $ and $ \{(x,y)~|~ x=0,-0.35 \leq y \leq -0.15\} $ and two outlets at the same vertical positions along $x=1.5$. The flow at the inlet is specified by $\boldsymbol{g}=\left(-100\left(y^2-0.35^2\right)\left(y^2-0.15^2\right), 0\right)^{\mathrm{T}}$. \label{exp:bypass}

  \item (\textit{Diffuser in 3D}) The domain $\Omega$ is a cube $(0,1)^{3}$. The inlet and outlet boundary conditions are imposed on $x=0$ (with unit flux in the positive $x$-direction) and on $\{\left(x,y\right)~|~x=1,(y-0.5)^2+(z-0.5)^2\leq 0.04\}$, respectively. \label{exp:diffuser3d}

  \item  (\textit{Pipe in 3D}) The domain $\Omega$ is a cube $(0,1)^3$. Two inlets represent circular regions at both ends of $\Omega$ in the $x$-direction, positioned along the $x=0$ and $x=1$ planes, respectively. Each of these regions has a radius of 0.2, centered at the point $(y=0.5, z=0.7)$ in the $y z$-plane with $\bold{g}=\left(\pm 1,0,0\right)^{\mathrm{T}}$. Similarly, the remaining two inlets are located at the boundaries in the $y$-direction, positioned along the $y=0$ and $y=1$ planes, with their centers at $(x=0.5, z=0.7)$ in the $x z$-plane given $\bold{g}=\left(0,\pm 1,0\right)^{\mathrm{T}}$. The outlet is defined as a circular region in the $x y$-plane at $z=0$ with the same radius, centered at the point $(x=0.5, y=0.5)$. \label{exp:pipe3d}
\end{enumerate}

\begin{table}[hbt!]
\centering
\begin{threeparttable}
\caption{The hyper-parameters for Algorithm \ref{Alg:UniformRefinement} and the phase--field model \eqref{eq:phasefieldupdateG1} used in the experiments. \label{tab:hyper}}
\begin{tabular}{c|ccccccc}
\toprule
Case & $(N,M)$ & $\Delta t$ & $\epsilon$  & $\gamma$  & $\zeta_{0}$ & $\beta$ & $\widetilde{S}$ \\
\midrule
 (a) & $(50,10)$ & 5e-4 & 1e-2  & 1e-2  & $100$ & $0.3$    & $0.25$\\
 (b) & $(50,10)$ & 1e-4 & 1e-2  & 1e-2  & $100$ & $0.5$    & $0.25$\\
 (c) & $(50,10)$ & 5e-5 & 1e-2  & 1e-2  & $100$ & $0.36$   & $0.25$\\
 (d) & $(50,10)$ & 1e-3 & 1e-3  & 1e-3  & $100$ & $0.925$  & $0.25$\\
 (e) & $(50,10)$ & 5e-4 & 5e-3  & 1e-1  & $50$  & $0.1667$ & $1$\\
 (f) & $(25,10)$ & 5e-4 & 1e-2  & 1e-2  & $500$ & $0.3$    & $0.25$\\
 (g) & $(25,10)$ & 5e-3 & 1e-2  & 1e-2  & $500$ & $0.3$    & $0.25$\\
\bottomrule
\end{tabular}
\end{threeparttable}
\end{table}

\begin{figure}[htb!]
	\centering \setlength{\tabcolsep}{0pt}
	\begin{tabular}{cccc}
        \subfigure[cases \ref{exp:pipebend} and \ref{exp:threeinflows}]{\includegraphics[width=1.5in]{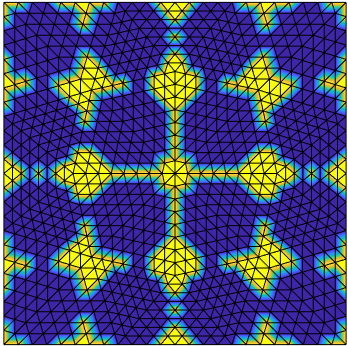}}
        &\subfigure[case \ref{exp:leftinflow}]{\includegraphics[width=1.5in]{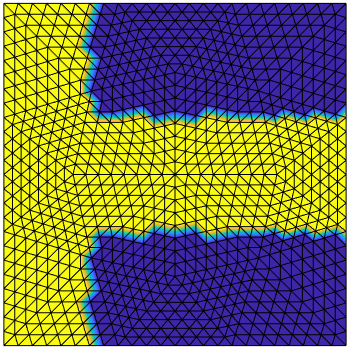}}
		&\subfigure[case \ref{exp:rugby}]{\includegraphics[width=2in]{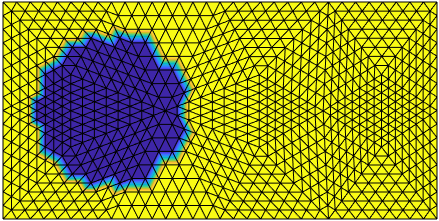}}
        &\subfigure[case \ref{exp:bypass}]{\includegraphics[width=2in]{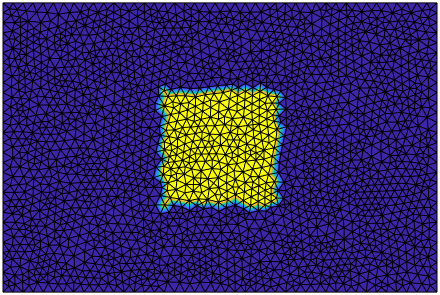}}
	\end{tabular}
\caption{The initial phase-field functions for the 2d cases.}
\label{fig: initialshape2d}
\end{figure}

The values of hyper-parameters (cf. Section \ref{subsect:solver}) are listed in Table \ref{tab:hyper}. In the table, $N$ and $M$ denote the number of iterations for solving the Stokes equation \eqref{min_phase-field}-\eqref{vp_stokes_div-free} and phase-field equation \eqref{eq:phasefieldupdateG1}, respectively. The parameter $\Delta t$ is the pseudo-time step size for updating the phase field. The three parameters ($\eps$, $\gamma$, $\widetilde{S}$) appear in the phase-field model \eqref{eq:phasefieldupdateG1}. The scalar $\beta$ represents the volume fraction. Throughout, the maximum number of uniform refinements $K$ is set to $3$, the penalty coefficient $\alpha_{0}$ is equal to $10000$, and the penalty scaling factor $\kappa$ is $1.1$. The Lagrange multiplier $\ell_{0}$ is initialized to zero. We compare the numerical results by CR-$P_{0}$ elements with those by $P_{2}$-$P_{1}$ Taylor-Hood elements. The initial mesh is indicated by $k=0$ and then uniformly refined three times. Table \ref{tab:mesh} lists the data of the meshes for the 2D cases, including the number of vertices, elements and degrees of freedom. During the refinement process, the CR-$P_{0}$ elements have fewer degrees of freedom than the $P_{2}$-$P_{1}$ elements. The initial phase-field functions in the 2d cases are shown in Figure \ref{fig: initialshape2d}, where yellow and blue represent values of $1$ and $0$, respectively.

\begin{table}[htb!]
\centering
\caption{Mesh data at different levels for the 2d cases.}\label{tab:mesh}
\begin{tabular}{c|c|cccccc}
\toprule
\multirow{2}{*}{Case} &
\multirow{2}{*}{mesh level} & \multirow{2}{*}{vertices} & \multirow{2}{*}{elements} & \multicolumn{2}{c}{dofs of CR-$P_0$} & \multicolumn{2}{c}{dofs of $P_2$-$P_1$}        \\ \cline{5-8}
 && & & $\bold{u}$ & $p$ & $\bold{u}$ & $p$ \\
 \midrule
\ref{exp:pipebend}-\ref{exp:threeinflows}&0 & 945 & 1792 & 5472 & 1792 & 7362 & 945 \\
&1 & 3681 & 7168 & 21696 & 7168 & 29058 & 3681 \\
&2 & 14529 & 28672 & 86400 & 28672 & 115458 & 14529 \\
&3 & 57729 & 114688 & 344832 & 114688 & 460290 & 57729 \\
\midrule
\ref{exp:rugby} & 0 & 881 & 1664 & 5088 & 1664 & 6850 & 881 \\
&1 & 3425 & 6656 & 20160 & 6656 & 27010 & 3425 \\
&2 & 13505 & 26624 & 80256 & 26624 & 107266 & 13505 \\
&3 & 53633 & 106496 & 320256 & 106496 & 427522 & 53633 \\
\midrule
\ref{exp:bypass}&0 & 2174 & 4202 & 12750 & 4202 & 17098 & 2174 \\
&1 & 8549 & 16808 & 50712 & 16808 & 67810 & 8549 \\
&2 & 33905 & 67232 & 202272 & 67232 & 270082 & 33905 \\
&3 & 53633 & 268928 & 807936 & 268928 & 1078018 & 135041 \\
\bottomrule
\end{tabular}

\end{table}


\begin{figure}[htb!]
	\centering\setlength{\tabcolsep}{0pt}
	\begin{tabular}{cccc}
        \includegraphics[width=1.5in]{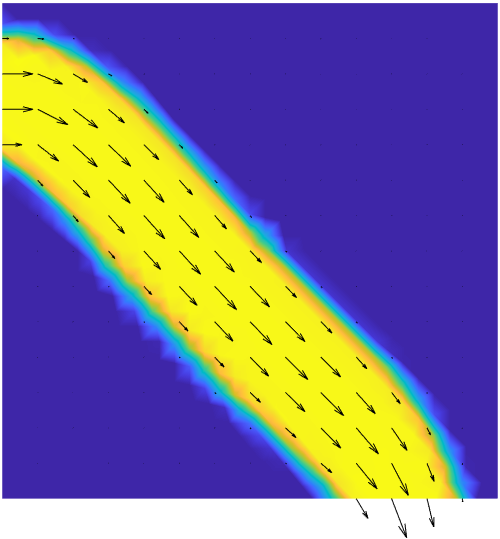}
        &\includegraphics[width=1.5in]{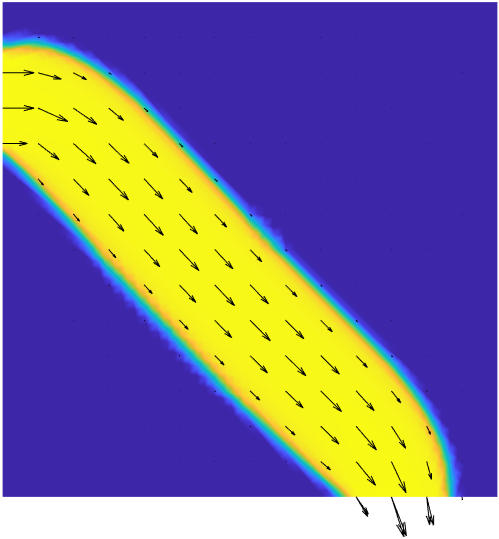}
        &\includegraphics[width=1.5in]{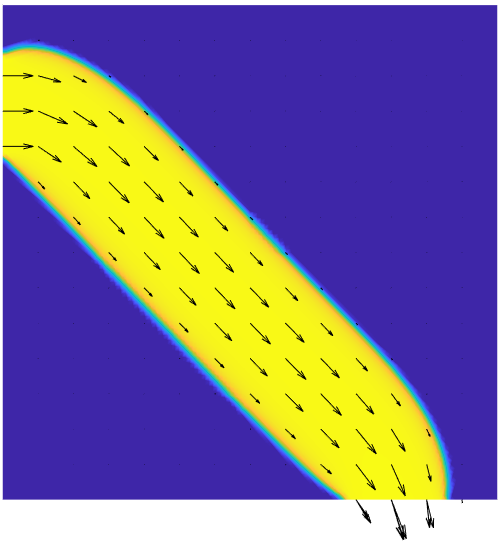}
        &\includegraphics[width=1.5in]{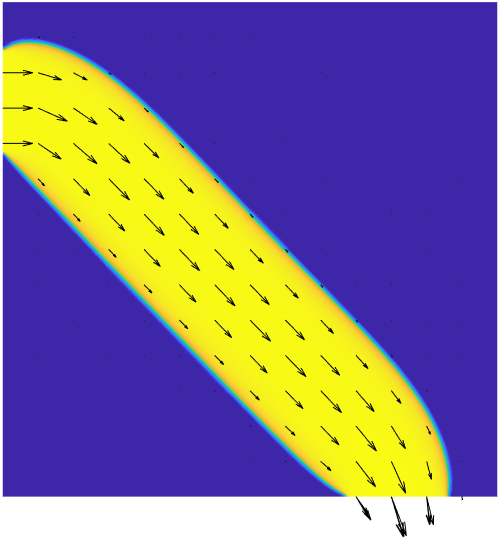}
        \\
        \includegraphics[width=1.5in]{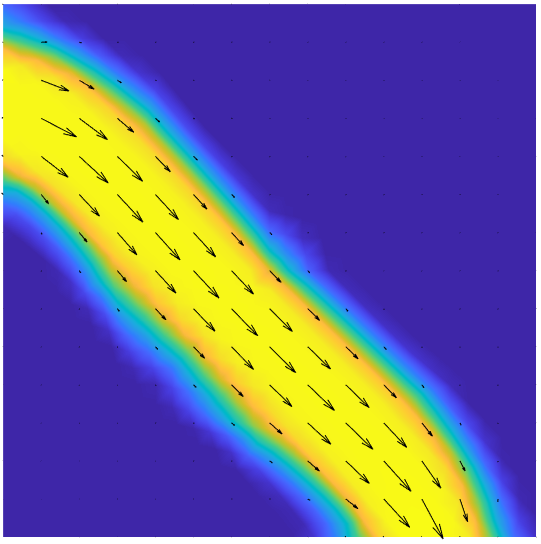}
        &\includegraphics[width=1.5in]{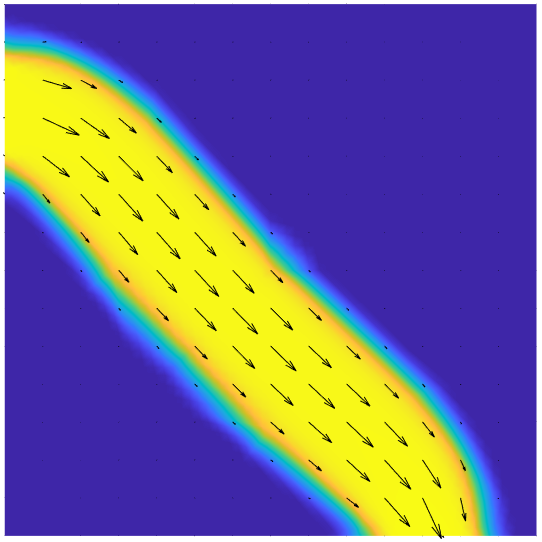}
        &\includegraphics[width=1.5in]{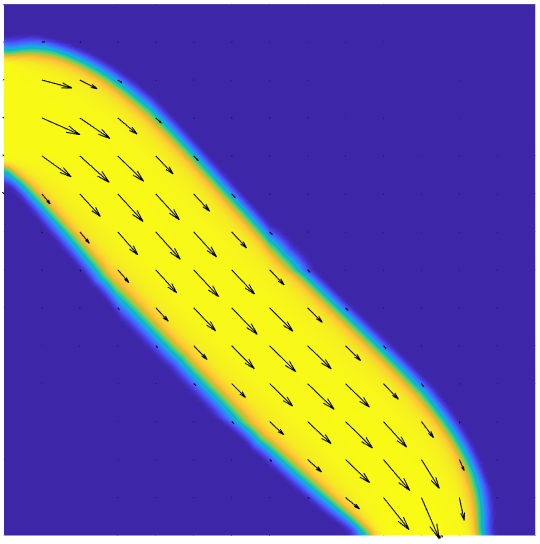}
        &\includegraphics[width=1.5in]{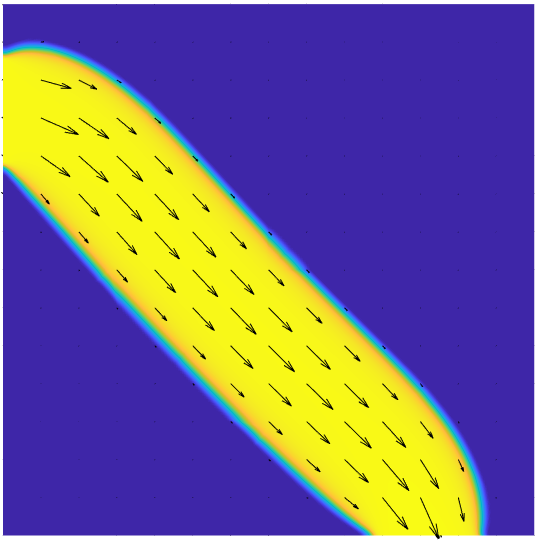}
	\end{tabular}
\caption{The optimal designs $\phi_k^\ast$ over each $\cT_k$ for case \ref{exp:pipebend} by CR-$P_0$ (1st row) and $P_2$-$P_1$ (2nd row).}
\label{fig:pipebendoptimizationresults}
\end{figure}

\begin{figure}[htb!]
	\centering\setlength{\tabcolsep}{0pt}
	\begin{tabular}{cccc}
        \includegraphics[width=1.5in]{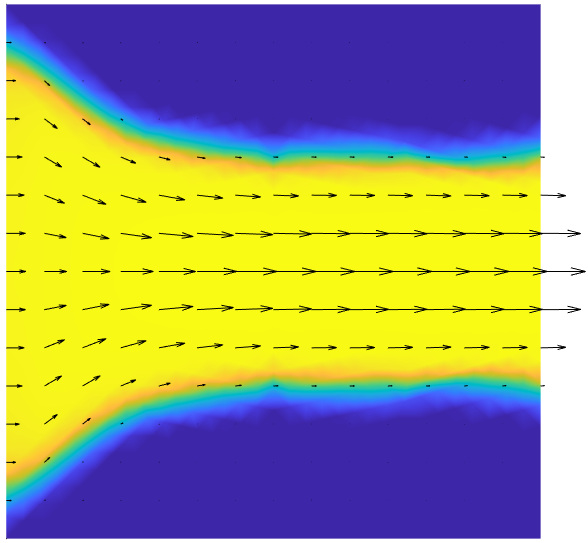}
        &\includegraphics[width=1.5in]{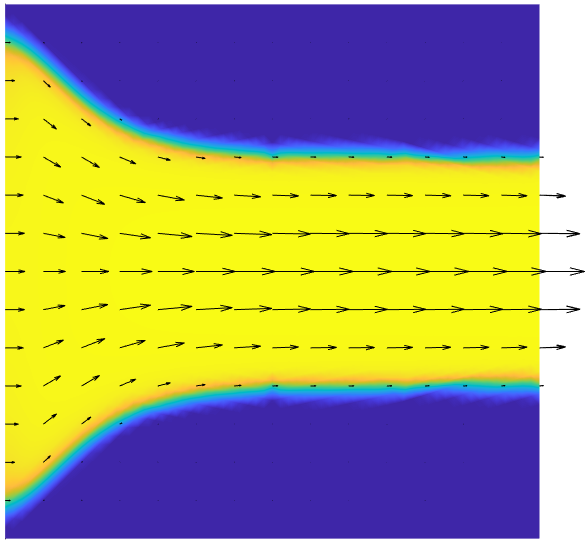}
        &\includegraphics[width=1.5in]{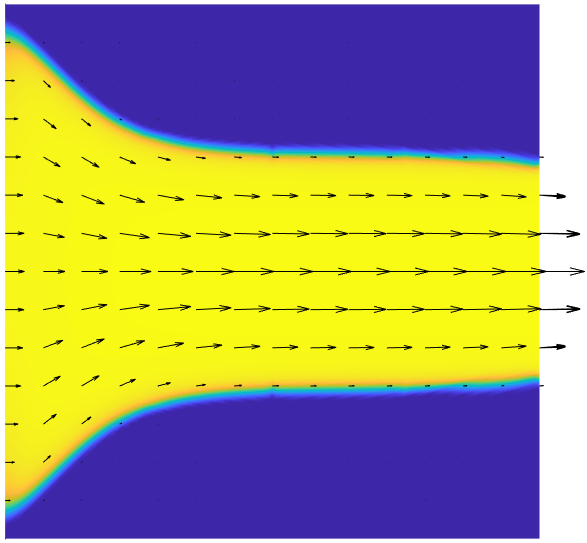}
        &\includegraphics[width=1.5in]{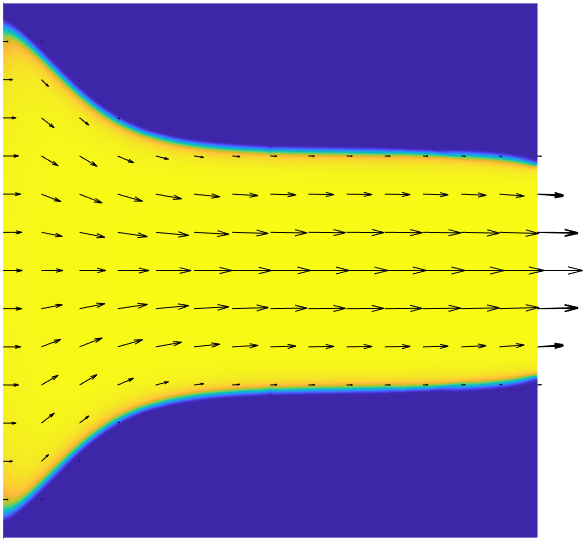}
        \\
        \includegraphics[width=1.5in]{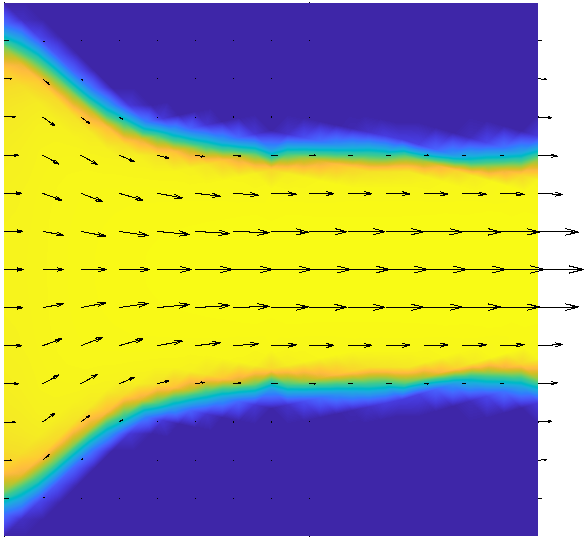}
        &\includegraphics[width=1.5in]{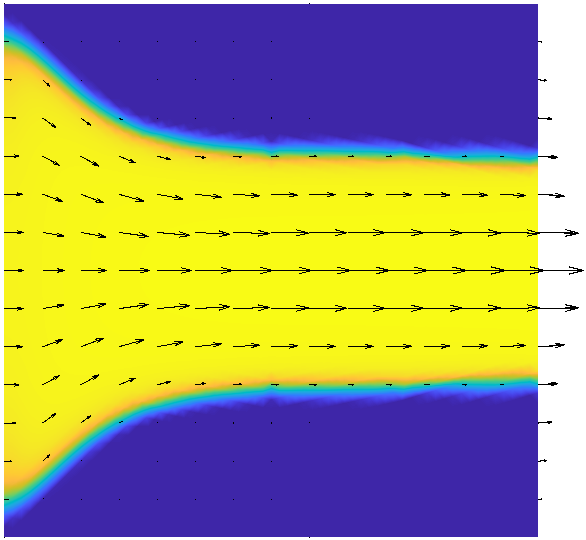}
        &\includegraphics[width=1.5in]{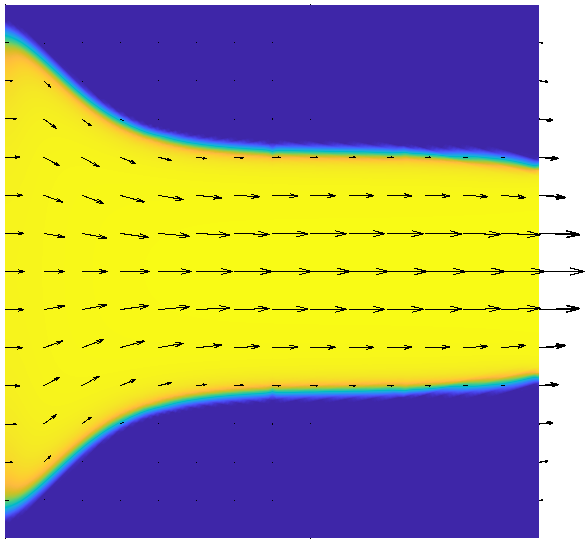}
        &\includegraphics[width=1.5in]{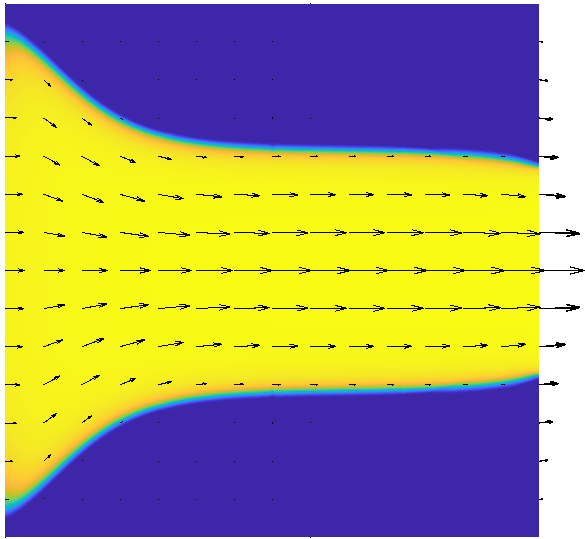}
	\end{tabular}
\caption{The optimal designs $\phi_k^\ast$  over each $\cT_k$ for case \ref{exp:leftinflow} by CR-$P_0$ (1st row) and  $P_2$-$P_1$ (2nd row).}
\label{fig:leftinflowoptimizationresults}
\end{figure}

\begin{figure}[htb!]
	\centering\setlength{\tabcolsep}{0pt}
	\begin{tabular}{cccc}
        \includegraphics[width=1.5in]{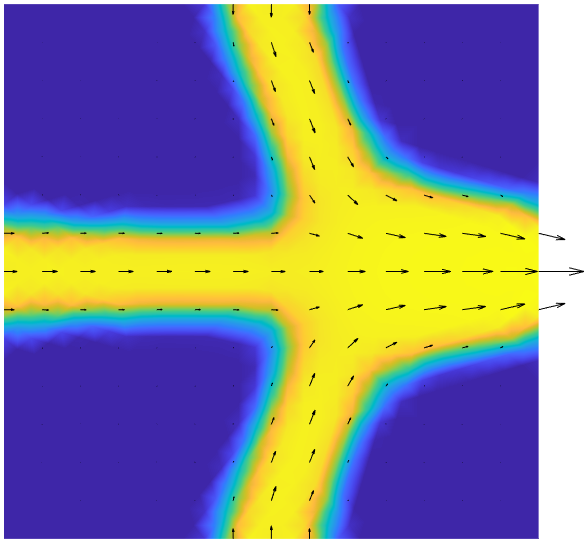}
        &\includegraphics[width=1.5in]{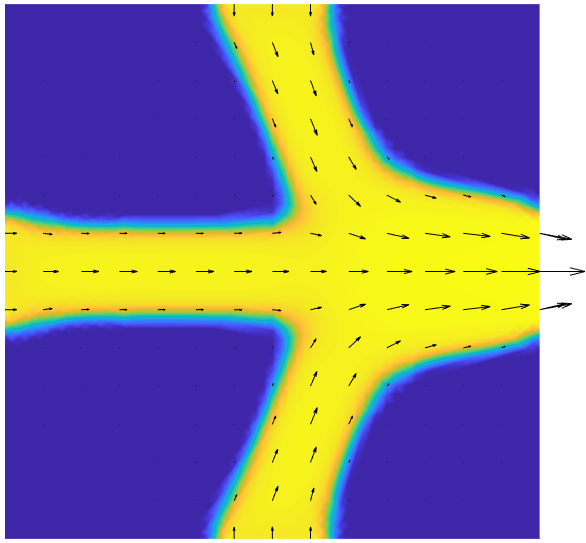}
        &\includegraphics[width=1.5in]{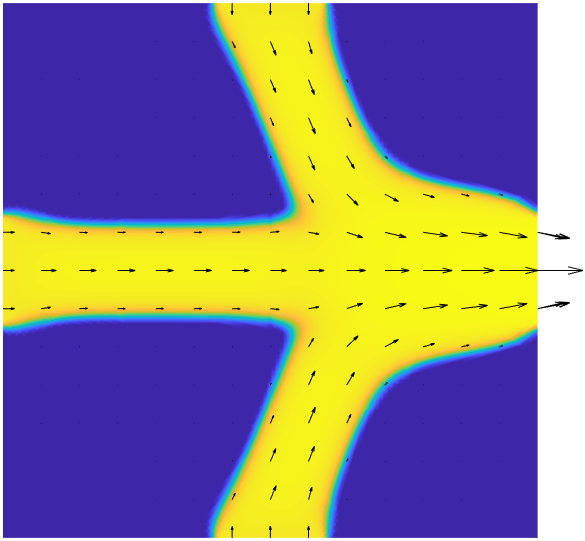}
        &\includegraphics[width=1.5in]{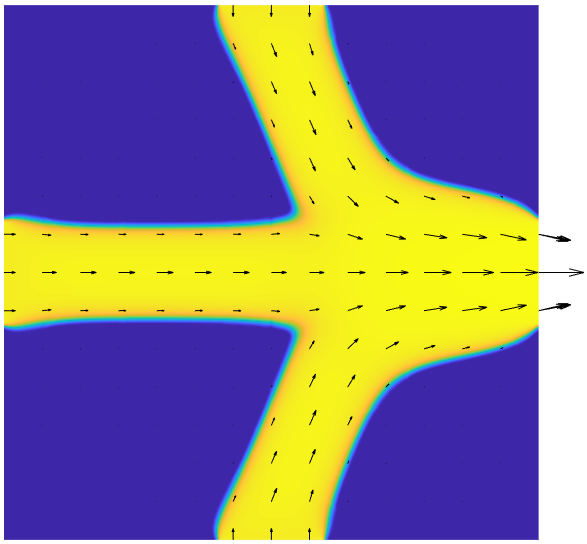}
        \\
        \includegraphics[width=1.5in]{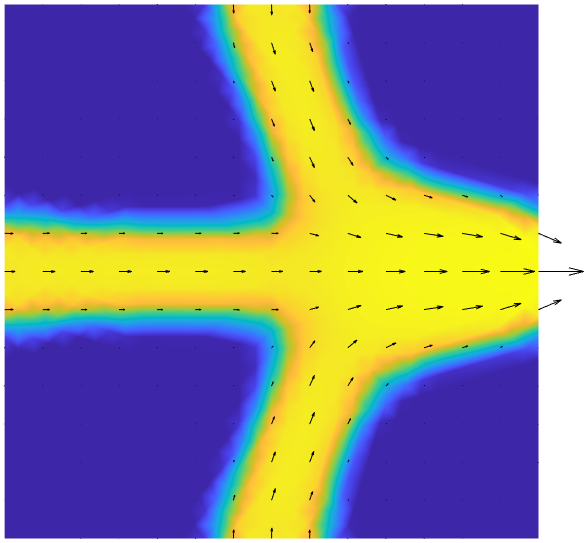}
        &\includegraphics[width=1.5in]{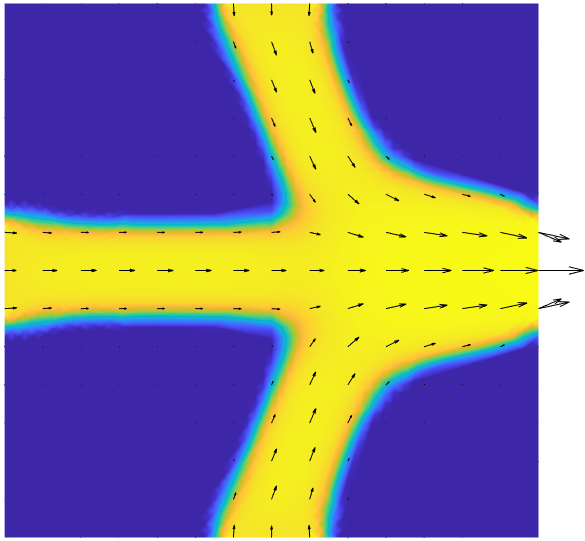}
        &\includegraphics[width=1.5in]{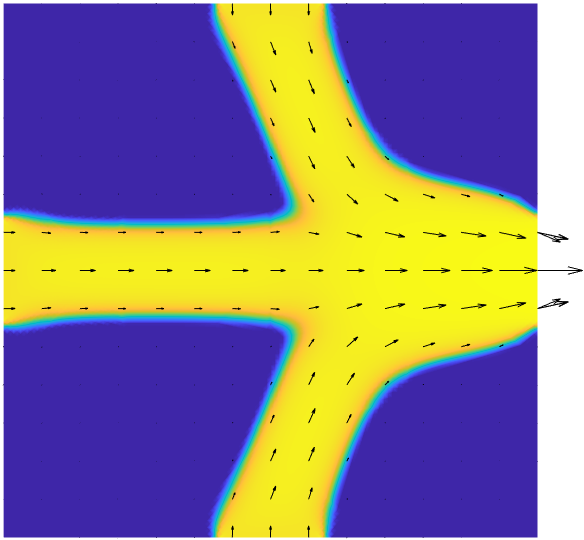}
        &\includegraphics[width=1.5in]{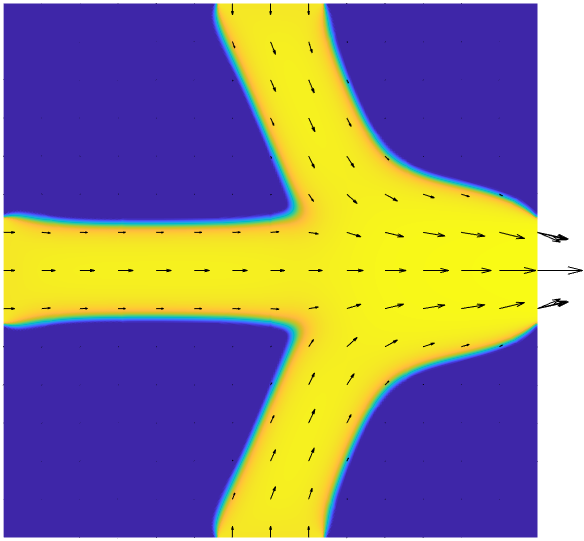}
	\end{tabular}
\caption{The optimal designs $\phi_k^\ast$ over each $\cT_k$ for case  \ref{exp:threeinflows} by CR-$P_0$ (1st row) and $P_2$-$P_1$ (2nd row).}
\label{fig:threeinflowsoptimizationresults}
\end{figure}

\begin{figure}[htb!]
	\centering\setlength{\tabcolsep}{0pt}
	\begin{tabular}{cccc}
        \includegraphics[width=1.5in]{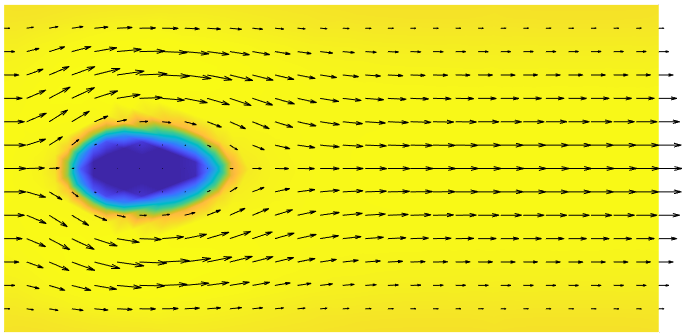}
        &\includegraphics[width=1.5in]{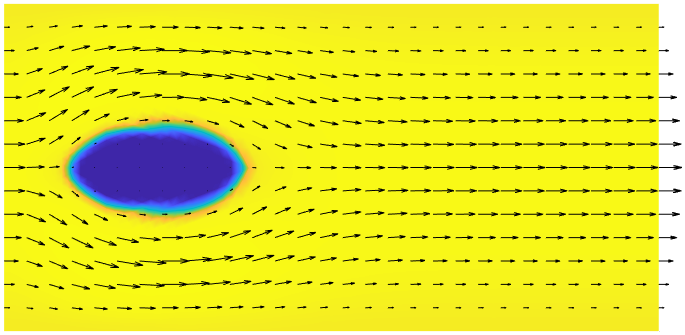}
        &\includegraphics[width=1.5in]{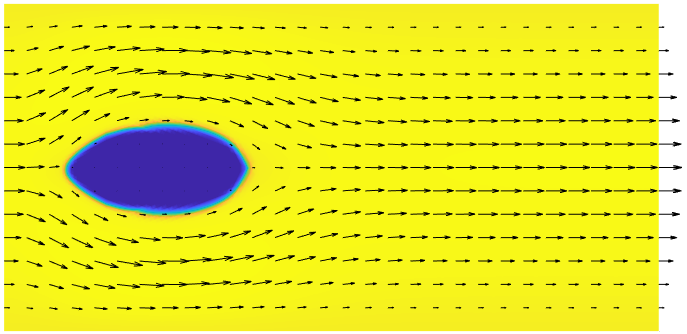}
        &\includegraphics[width=1.5in]{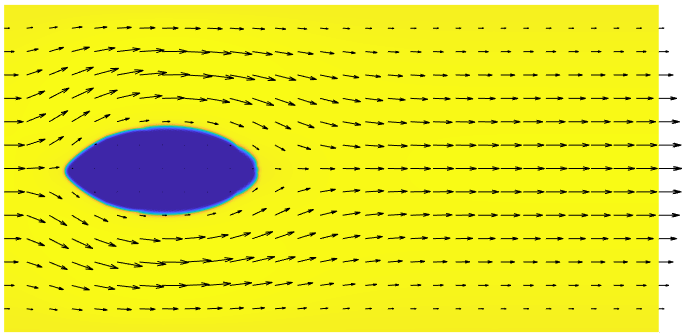}
        \\
        \includegraphics[width=1.5in]{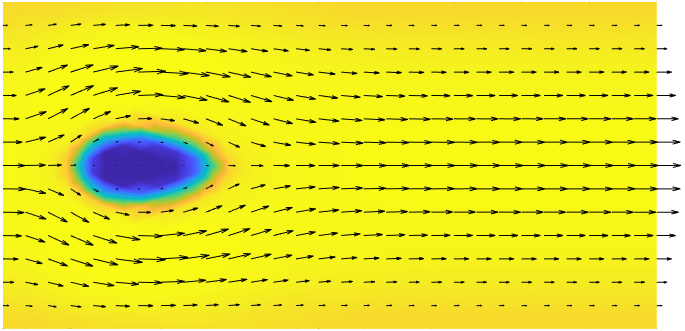}
        &\includegraphics[width=1.5in]{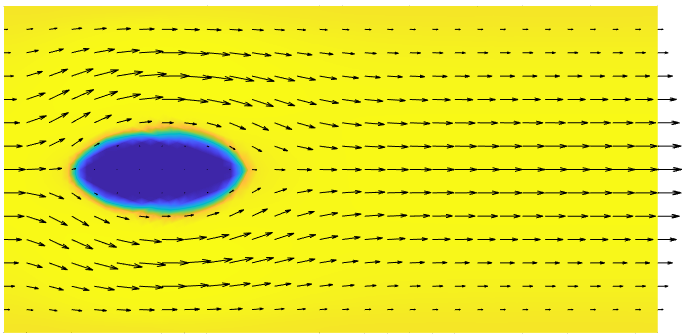}
        &\includegraphics[width=1.5in]{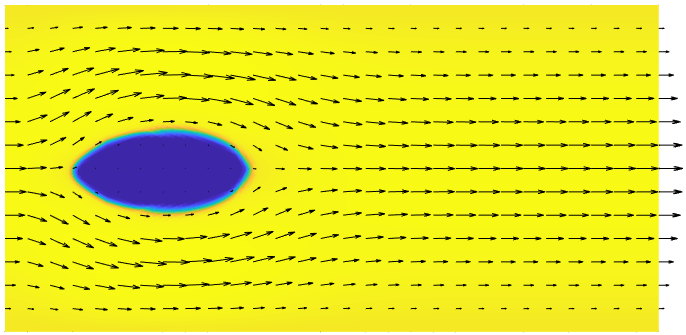}
        &\includegraphics[width=1.5in]{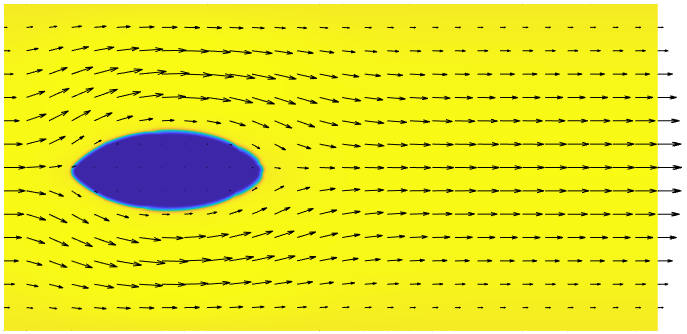}
	\end{tabular}
\caption{The optimal designs $\phi_k^\ast$ over each $\cT_k$ for case \ref{exp:rugby} by CR-$P_0$ (1st row) and $P_2$-$P_1$ (2nd row).}
\label{fig:rugbyoptimizationresults}
\end{figure}

\begin{figure}[htb!]
\centering\setlength{\tabcolsep}{0pt}
	\begin{tabular}{cccc}
        \includegraphics[width=1.5in]{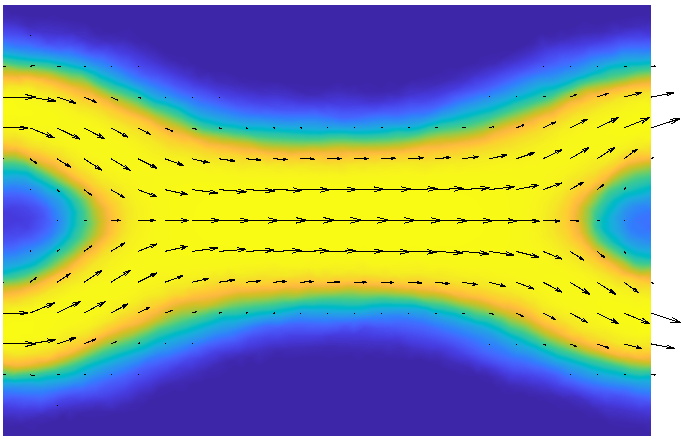}
        &\includegraphics[width=1.5in]{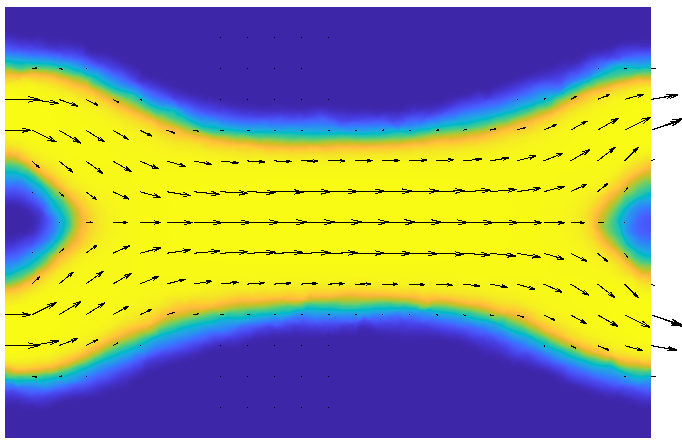}
        &\includegraphics[width=1.5in]{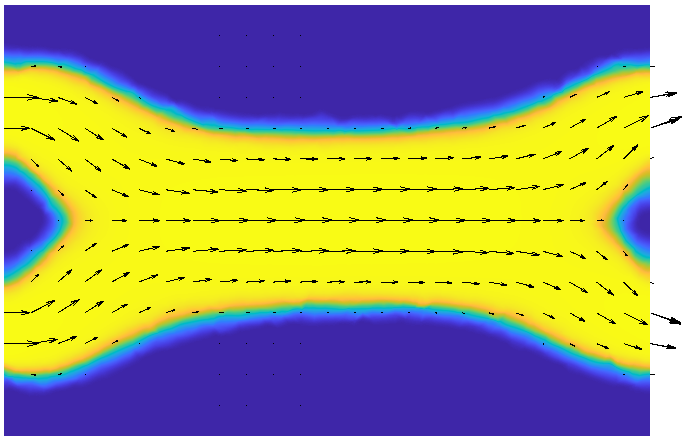}
        &\includegraphics[width=1.5in]{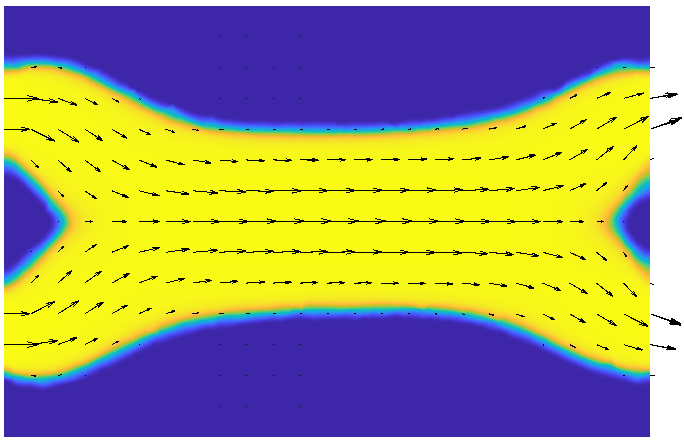}
        \\
        \includegraphics[width=1.5in]{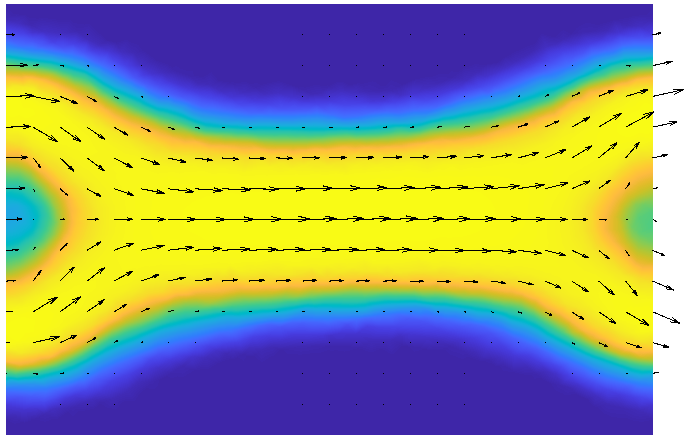}
        &\includegraphics[width=1.5in]{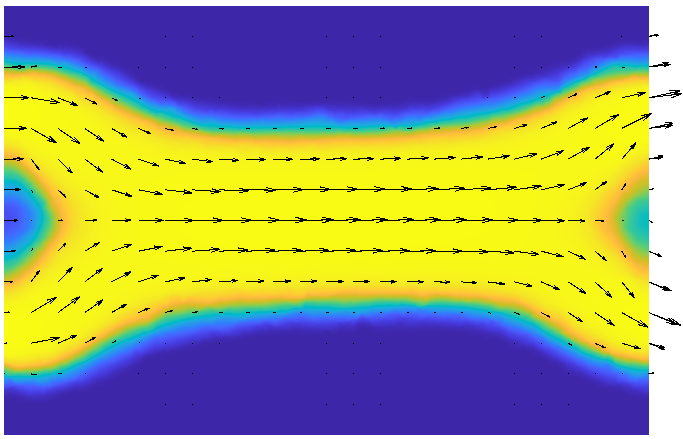}
        &\includegraphics[width=1.5in]{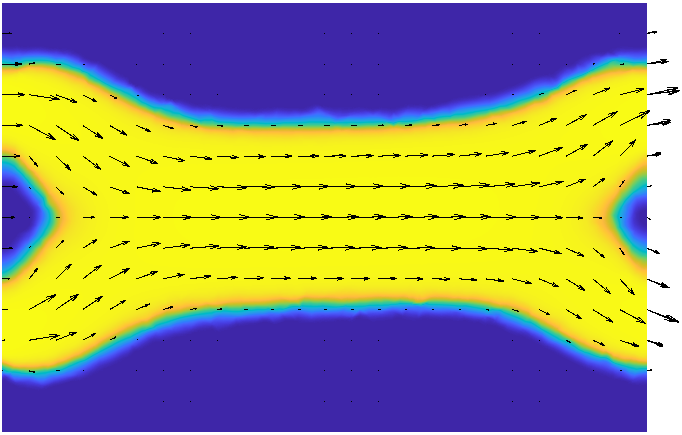}
        &\includegraphics[width=1.5in]{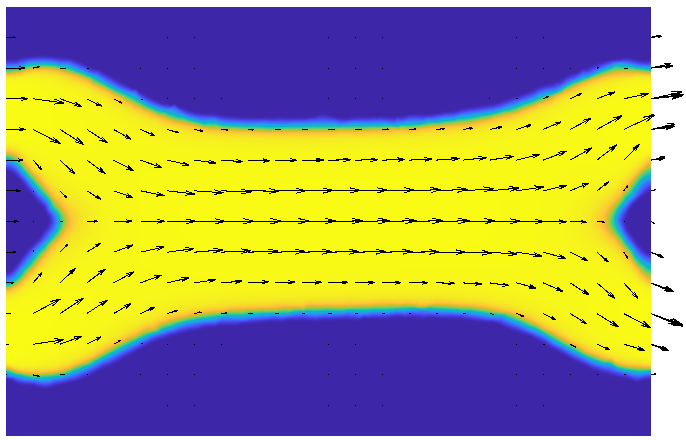}
	\end{tabular}
\caption{The optimal designs $\phi_k^\ast$ over each $\cT_k$ for case \ref{exp:bypass} by CR-$P_0$ (1st row) and $P_2$-$P_1$ (2nd row).}
\label{fig:bypassoptimizationresults}
\end{figure}

\begin{figure}[htb!]
\centering\setlength{\tabcolsep}{0pt}
\begin{tabular}{cc|cc}
\toprule
      \multicolumn{2}{c}{ \ref{exp:pipebend}} & \multicolumn{2}{c}{\ref{exp:leftinflow}}\\
        \includegraphics[width=1.6in]{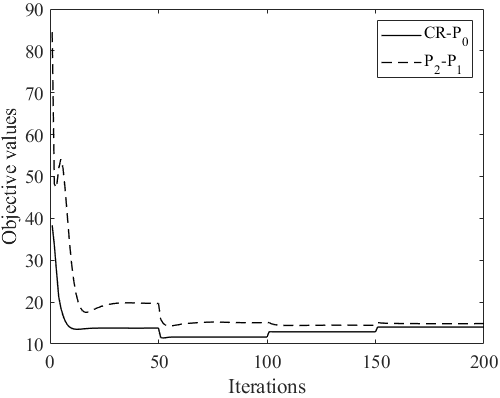}
        & \includegraphics[width=1.6in]{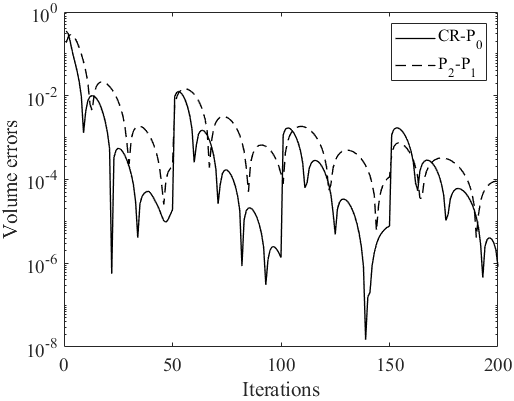} &\includegraphics[width=1.6in]{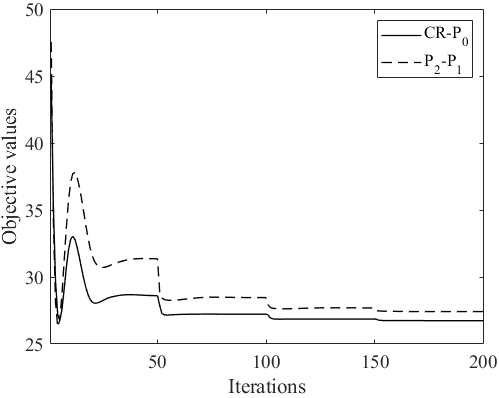} &\includegraphics[width=1.6in]{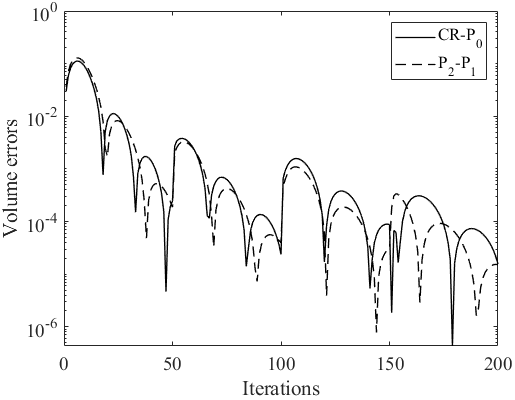}\\
    \midrule
   \multicolumn{2}{c}{\ref{exp:threeinflows}} & \multicolumn{2}{c}{\ref{exp:rugby}}\\ \includegraphics[width=1.6in]{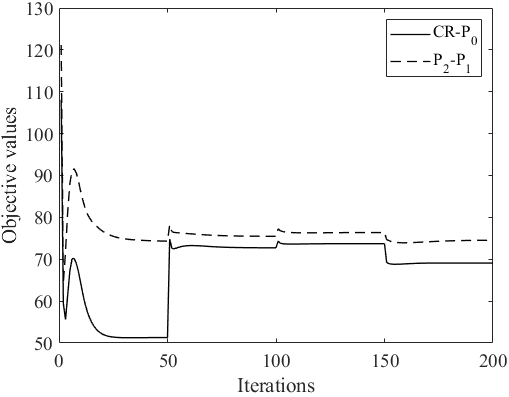}
        &\includegraphics[width=1.6in]{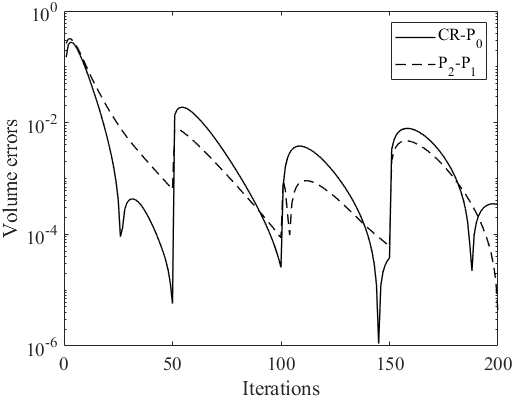} & \includegraphics[width=1.6in]{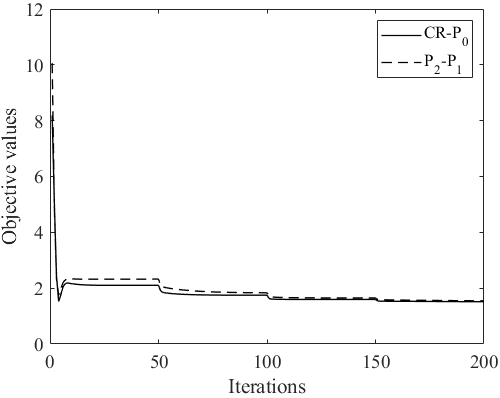}
    &\includegraphics[width=1.6in]{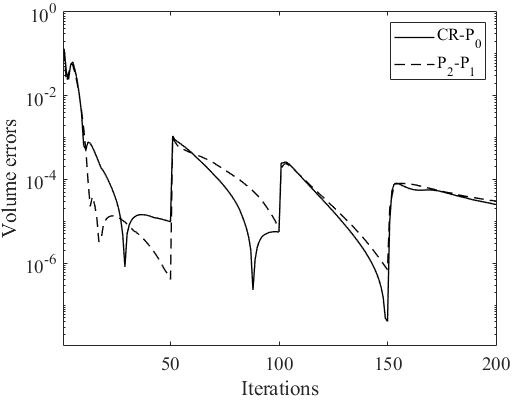}\\
    \midrule
    \multicolumn{2}{c}{\ref{exp:bypass}} & \multicolumn{2}{c}{\ref{exp:diffuser3d}}\\
\includegraphics[width=1.6in]{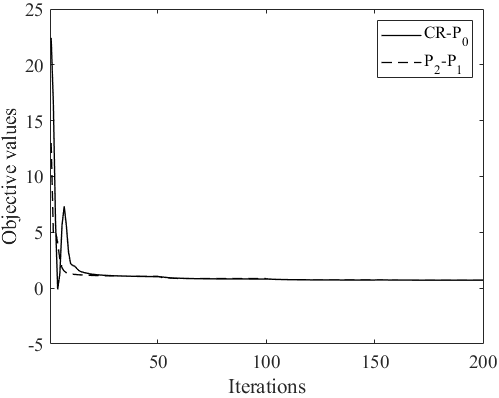}
&\includegraphics[width=1.6in]{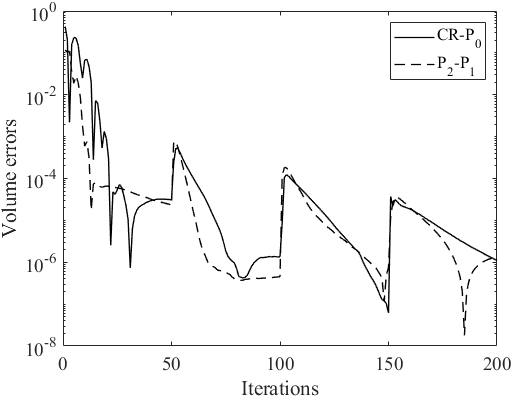} &
\includegraphics[width=1.6in]{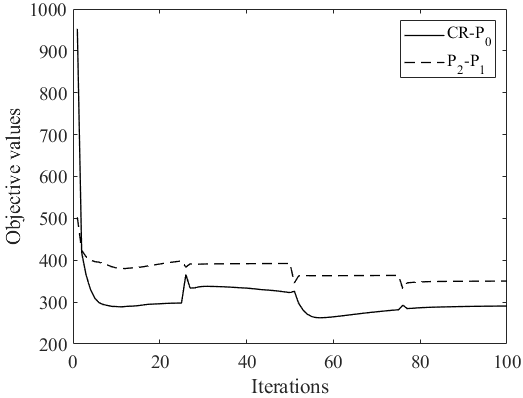} &\includegraphics[width=1.6in]{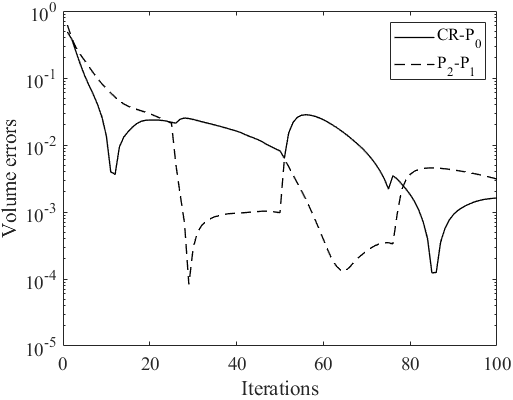}\\
\midrule
\multicolumn{2}{c}{\ref{exp:pipe3d}}\\
\includegraphics[width=1.6in]{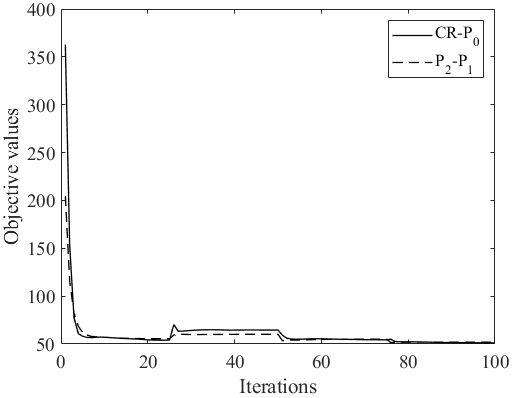}
&\includegraphics[width=1.6in]{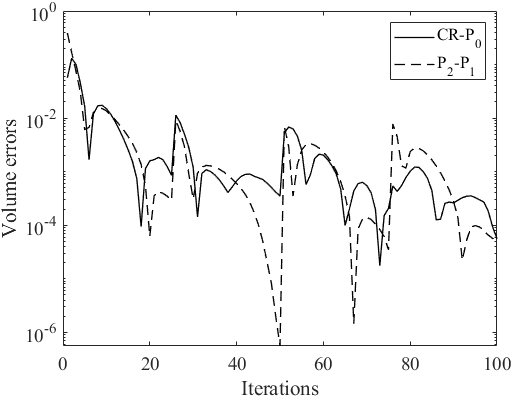}\\
\bottomrule
\end{tabular}
\caption{The convergence history of the objective value and the volume constraint error as functions of the total number ($K \times N$) of outer-iterations.\label{fig:evol-obj}}
\end{figure}

Figures \ref{fig:pipebendoptimizationresults}-\ref{fig:bypassoptimizationresults} display the optimal designs for the 2d cases, which are presented from left to right as the refinement proceeds. The first and second rows are due to the CR-$P_{0}$ elements and $P_2$-$P_{1}$ Taylor-Hood elements, respectively. The interfaces between the fluid and void regions progressively sharpens as the refinement proceeds: the optimization on the initial mesh ($k=0$) yields a rough shape of the optimal design, and the subsequent refinements ($k=1,2,3$) sharpen the interfaces of the fluid, due to the better resolution of the velocity field $\bold{u}$ by both CR-$P_{0}$ and $P_{2}$-$P_{1}$ elements. This shows the effectiveness of Algorithm \ref{Alg:UniformRefinement}. The convergence history of the total energy and the volume constraint error is illustrated in Figure \ref{fig:evol-obj}. The horizontal axis represents the total iterations on four mesh levels, and we observe the decay of the objective value $\mathcal{L}(\phi,\bold{u})$ and the well-preserved volume error by CR-$P_{0}$ and $P_{2}$-$P_{1}$ elements. Also, the total energy decreases rapidly at the early stage of iterations, and then stabilizes. The interpolation and projection on different mesh levels cause slight spikes with the total energy in each case. 

\begin{figure}[htbp!]
	\centering\setlength{\tabcolsep}{0pt}
	\begin{tabular}{cccc}
		\subfigure[case \ref{exp:diffuser3d}]{\includegraphics[width=1.5in]{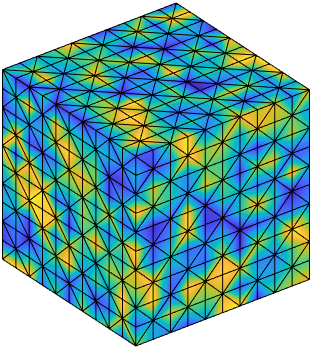}}
        &\subfigure[case \ref{exp:diffuser3d}]{\includegraphics[width=1.5in]{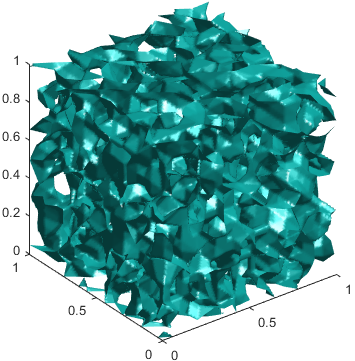}}
        \hspace{2mm}
		&\subfigure[case \ref{exp:pipe3d}]{\includegraphics[width=1.5in]{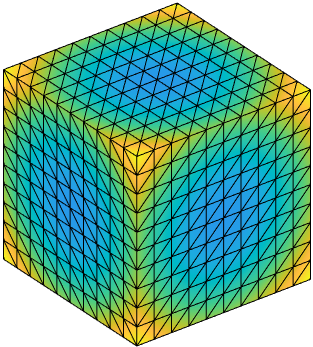}}
        &\subfigure[case \ref{exp:pipe3d}]{\includegraphics[width=1.5in]{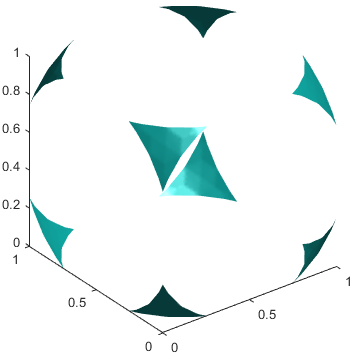}}
	\end{tabular}
\caption{The initial phase-field functions and designs for cases \ref{exp:diffuser3d} and \ref{exp:pipe3d}. Blue and yellow denote zero and one, respectively. (i) and (iii): the phase-field functions. (ii) and (iv): associated iso-surfaces with a value of $0.5$.}
\label{fig: initialshape3d}
\end{figure}

\begin{figure}[htb!]
	\centering\setlength{\tabcolsep}{0pt}
	\begin{tabular}{cccc}
        \includegraphics[width=1.2in]{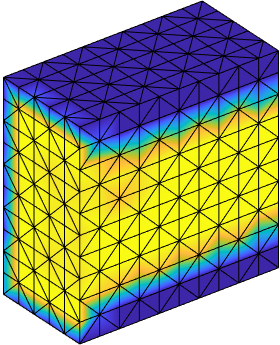}
        &\includegraphics[width=1.2in]{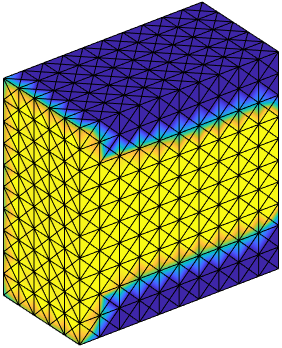}
        &\includegraphics[width=1.2in]{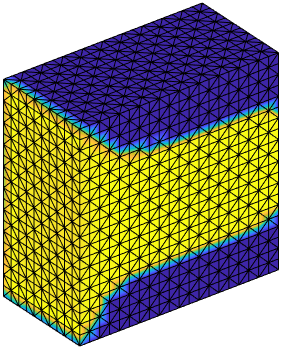}
        &\includegraphics[width=1.2in]{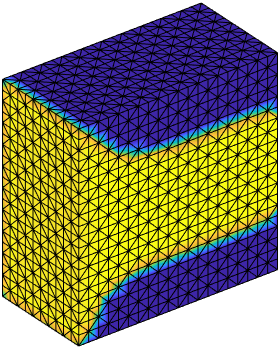}
        \\
        \includegraphics[width=1.7in]{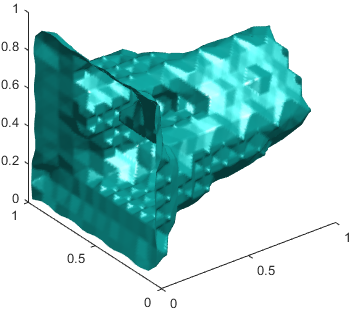}
        &\includegraphics[width=1.7in]{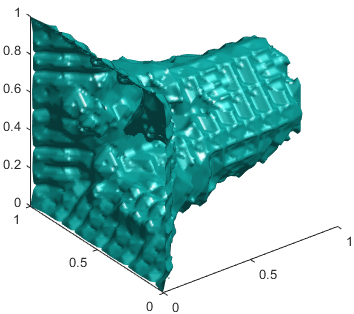}
        &\includegraphics[width=1.7in]{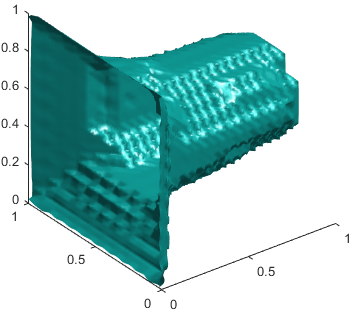}
        &\includegraphics[width=1.7in]{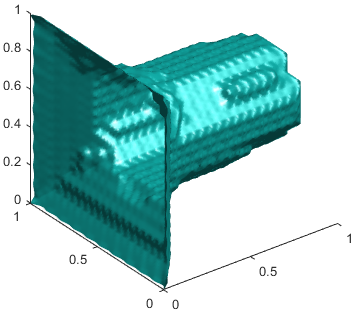}
        \\
        \includegraphics[width=1.2in]{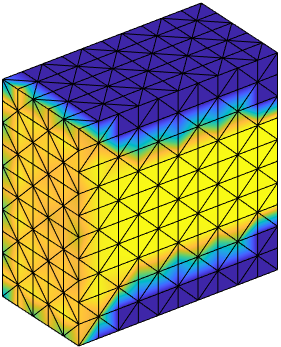}
        &\includegraphics[width=1.2in]{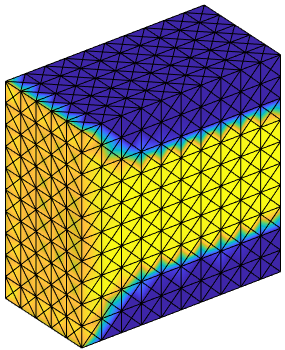}
        &\includegraphics[width=1.2in]{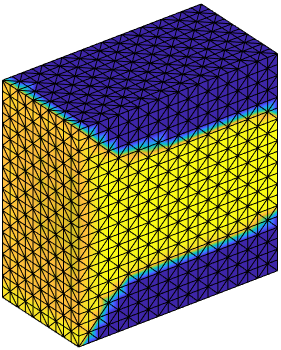}
        &\includegraphics[width=1.2in]{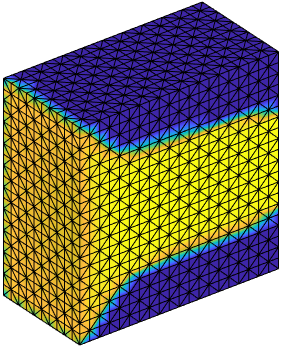}
        \\
        \includegraphics[width=1.7in]{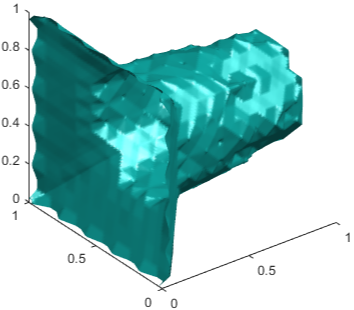}
        &\includegraphics[width=1.7in]{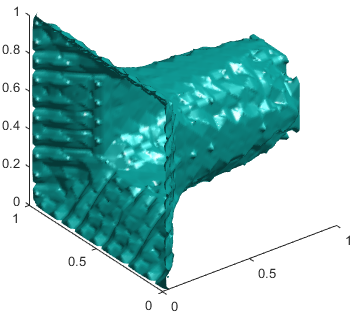}
        &\includegraphics[width=1.7in]{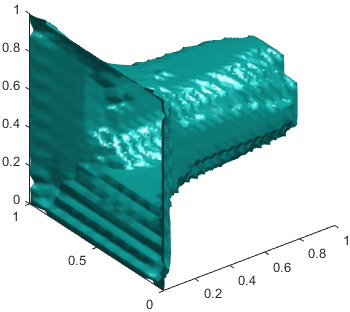}
        &\includegraphics[width=1.7in]{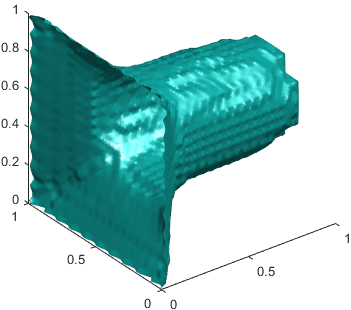}
	\end{tabular}
\caption{The optimal designs $\phi_k^\ast$ over each $\cT_k$ for case \ref{exp:diffuser3d} by CR-$P_0$ (first two rows) and $P_2$-$P_1$ (last two rows). The 1st and 3rd rows: phase-field functions. The 2nd and 4th rows: associated iso-surface with a value of 0.5.}
\label{fig: optimized results diffuser3d}
\end{figure}

\begin{figure}[htb!]
	\centering\setlength{\tabcolsep}{0pt}
	\begin{tabular}{cccc}
        \includegraphics[width=1.2in]{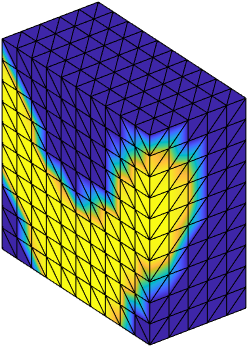}
        &\includegraphics[width=1.2in]{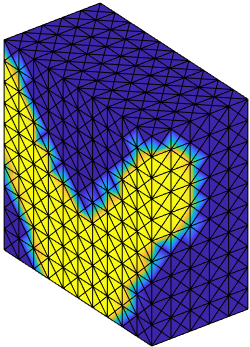}
        &\includegraphics[width=1.2in]{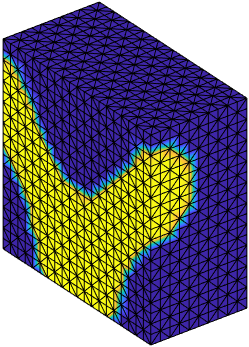}
        &\includegraphics[width=1.2in]{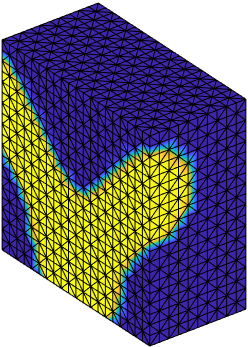}
        \\
        \includegraphics[width=1.7in]{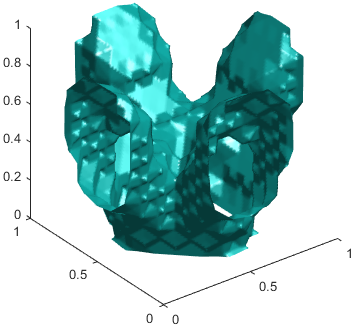}
        &\includegraphics[width=1.7in]{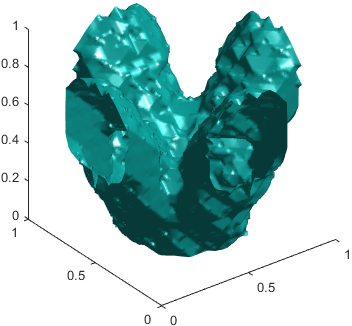}
        &\includegraphics[width=1.7in]{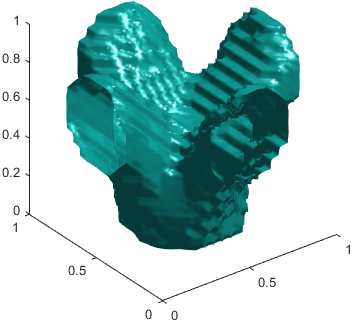}
        &\includegraphics[width=1.7in]{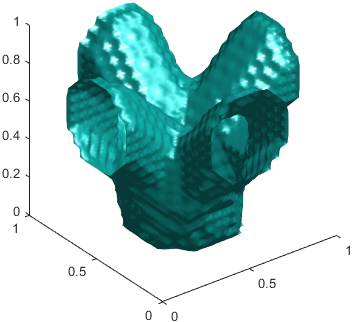}
        \\
        \includegraphics[width=1.2in]{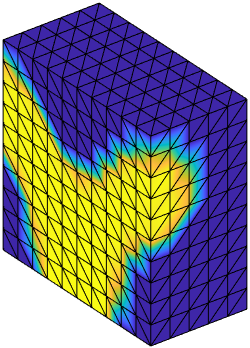}
        &\includegraphics[width=1.2in]{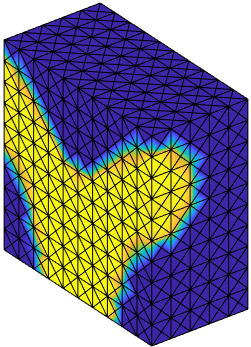}
        &\includegraphics[width=1.2in]{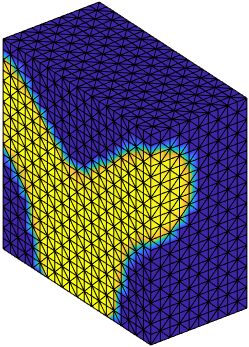}
        &\includegraphics[width=1.2in]{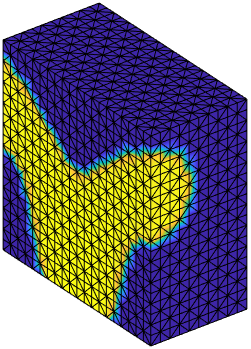}
        \\
        \includegraphics[width=1.7in]{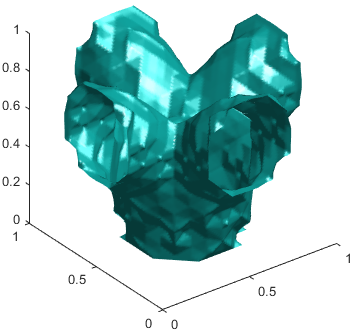}
        &\includegraphics[width=1.7in]{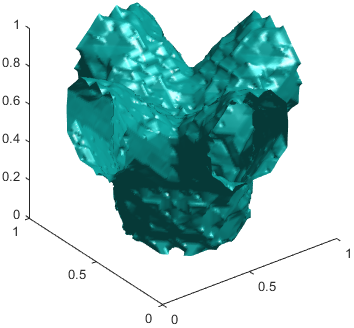}
        &\includegraphics[width=1.7in]{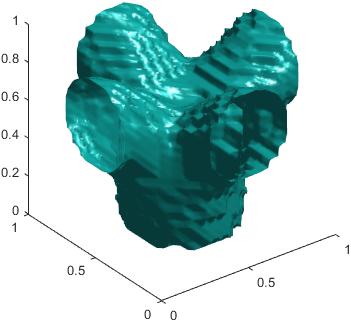}
        &\includegraphics[width=1.7in]{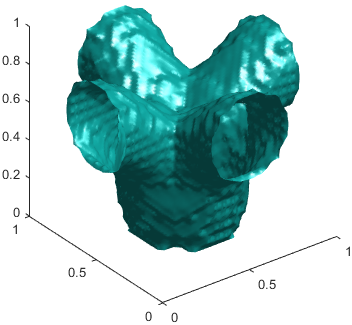}
	\end{tabular}
\caption{The optimal designs $\phi_k^\ast$ over each $\cT_k$ for case \ref{exp:pipe3d} by CR-$P_0$ (first two rows) and $P_2$-$P_1$ (last two rows). The 1st and 3rd rows: phase-field functions. The 2nd and 4th rows: associated iso-surface with a value of 0.5.}
\label{fig: optimized results pipe3d}
\end{figure}

The initial phase-field functions and topologies in the 3d cases are depicted in Figure \ref{fig: initialshape3d}, where the initial designs denote the associated iso-surface with a value of $0.5$. The initial phase-field is randomly generated with a uniform distribution in the interval $[0, 1]$ for case \ref{exp:diffuser3d}, and set to $(x-0.5)^2+(y-0.5)^2+(z-0.5)^2 - 0.25^2$ for case \ref{exp:pipe3d}. Table \ref{tab: diffuser pipe 3d data} presents the detailed information of mesh levels and objective values by CR-$P_{0}$ and $P_{2}$-$P_{1}$ elements in cases \ref{exp:diffuser3d}-\ref{exp:pipe3d}. In Figures \ref{fig: optimized results diffuser3d} and \ref{fig: optimized results pipe3d}, we summarize  optimal designs, including numerical phase-field functions and iso-surfaces with a value of $0.5$ on four mesh levels by CR-$P_{0}$ (first two rows) and $P_{2}$-$P_{1}$ (last two rows) elements. Similar to the 2d cases, the interfaces between fluid and void regions progressively sharpen as the mesh refinement proceeds. The convergence curves of the objective value and the volume error are presented in Figure \ref{fig:evol-obj}, indicating the good performance of the optimization process during mesh refinements.

\begin{table}[hbt!]
\centering
\caption{Mesh data and comparison of objective values by CR-$P_0$ and $P_2$-$P_1$ in cases \ref{exp:diffuser3d}-\ref{exp:pipe3d}.}
\label{tab: diffuser pipe 3d data}
\begin{tabular}{ccccccc}
\hline
\multirow{2}{*}{mesh level} & \multirow{2}{*}{vertices} & \multirow{2}{*}{elements} & \multicolumn{2}{c}{Case \ref{exp:diffuser3d}} & \multicolumn{2}{c}{Case \ref{exp:pipe3d}} \\ \cline{4-7}
 & & & obj (CR-$P_0$) & obj ($P_2$-$P_1$) & obj (CR-$P_0$) & obj ($P_2$-$P_1$) \\ \hline
0 & 2331 & 12000  & 297.79 & 398.17 & 54.07 & 55.59  \\
1 & 5631 & 24000  & 322.78 & 392.21 & 64.53 & 60.15  \\
2 & 9261 & 48000  & 281.82 & 363.78 & 54.42 & 55.07  \\
3 & 17261 & 96000 & 290.68 & 350.28 & 51.10 & 51.92  \\ \hline
\end{tabular}
\end{table}

\begin{table}[hbt!]
\centering
\caption{The comparison of the optimization results by CR-$P_0$ and $P_2$-$P_1$.}
\label{tab:comparisonallresults}
\begin{tabular}{cccccc}
\hline
\multirow{2}{*}{Case} & \multicolumn{2}{c}{CR-$P_0$} & \multicolumn{2}{c}{$P_2$-$P_1$} & \multirow{2}{*}{savings} \\ \cline{2-5}
                       & obj    & time (secs)   & obj    & time (secs)   &         \\ \hline
\ref{exp:pipebend}     & 14.04              & 494           & 14.87              & 1388          & 64.45\% \\
\ref{exp:leftinflow}   & 26.74              & 449           & 27.42              & 2036          & 77.95\% \\
\ref{exp:threeinflows} & 69.09              & 399           & 74.54              & 1084          & 63.20\% \\
\ref{exp:rugby}        & 1.51               & 424           & 1.55               & 1972          & 78.48\% \\
\ref{exp:bypass}       & 0.71               & 1029          & 0.72               & 2968          & 65.32\% \\
\ref{exp:diffuser3d}   &290.68              &1937           & 350.28             & 5064          & 61.74\% \\
\ref{exp:pipe3d}       &51.10               &2004           & 51.92              & 5273          & 61.98\% \\ \hline
\end{tabular}
\end{table}

Table \ref{tab:comparisonallresults} presents a detailed quantitative comparison of numerical results by CR-$P_0$ and $P_2$-$P_1$ elements for cases \ref{exp:pipebend}-\ref{exp:pipe3d}. The objective values produced by the CR-$P_0$ elements are very close to or slightly lower than that by the $P_2$-$P_1$ elements while the optimization process by the former saves over $60\%$ of the time compared to that by the latter, due to the fewer degrees of freedom of CR-$P_0$ on the same mesh level.

\paragraph{Acknowledgments} The work of B. Jin is supported by Hong Kong RGC General Research Fund (Project 14306824), and a start-up fund from The Chinese University of Hong Kong. The work of Y. Xu is in part supported by National Natural Science Foundation of China (12250013, 12261160361 and 12271367), Science and Technology Commission of Shanghai Municipality (20JC1413800 and 22ZR1445400) and General Research Grant (KF2023018 and KF2024068) from Shanghai Normal University. The work of S. Zhu is supported in part National Natural Science Foundation of China under grant 12471377 and Science and the Technology Commission of Shanghai Municipality (22DZ2229014).

\bibliographystyle{abbrv}
\bibliography{top_fluid}

\end{document}